\documentclass[12pt,notitlepage,twoside]{article}
\usepackage{latexsym}
\usepackage{amsmath}
\usepackage{amsfonts}
\usepackage{amssymb}
\renewcommand{\a }{\alpha }

\renewcommand{\d}{\delta }

\newcommand{\D }{\Delta }

\newcommand{\e }{\varepsilon }

 \newcommand{\foral }{\forall\, }

\renewcommand{\l }{\lambda }
\renewcommand{\L }{\Lambda }

\newcommand{\n }{\nabla }
\newcommand{\var }{\varphi }

\renewcommand{\th }{\theta }
\renewcommand{\o }{\omega }
\renewcommand{\O }{\Omega }

\newcommand{\be}{\begin{equation}}
\newcommand{\ee}{\end{equation}}
\newenvironment{pf}{\noindent{\bf Proof.}\enspace}{
\hfill$\Box$\medskip}
\newenvironment{pfn}[1]{\noindent{\bf Proof of {#1}\enspace}}{
\hfill$\Box$\medskip}
\newcommand{\R}{\mathbb{R}}
 
\newcommand{\N}{\mathbb{N}}

\newtheorem{thm}{Theorem}[section]
\newtheorem{pro}[thm]{Proposition}
\newtheorem{lem}[thm]{Lemma}
\newtheorem{rem}[thm]{Remark}

\numberwithin{equation}{section}

\author{Mohamed Ben Ayed$^a$\footnote{ Corresponding author. Fax : +216-74-274437, E-mail : \texttt{Mohamed.Benayed@fss.rnu.tn} },\, Khalil El Mehdi$^{b,c}$ \& Mokhless Hammami$^a$\,\footnote{E-mail adresses : \texttt{khalil@univ-nkc.mr} (K. El Mehdi), \texttt{Mokhless.Hammami@fss.rnu.tn} (M. Hammami).} \\
{\footnotesize
a : D{\'e}partement de Math{\'e}matiques, Facult{\'e} des Sciences de Sfax, Route
Soukra, Sfax, Tunisia.}\\
{\footnotesize
b :  Facult\'e des Sciences et Techniques, Universit\'e de Nouakchott, BP 5026, Nouakchott, Mauritania.}\\
{\footnotesize
 c : The Abdus Salam ICTP, Mathematics Section,
 Strada Costiera 11, 34014 Trieste, Italy.}
}

\title { \Large \textbf{ Nonexistence of Bounded Energy
Solutions for a Fourth Order  Equation on Thin Annuli }}

\begin{document}

\date{ }

\maketitle {\footnotesize \noindent {\bf ABSTRACT. }
 In this paper we study the problem  $
P_{\varepsilon}:\,
\Delta^2 u_\varepsilon =  u_\varepsilon^{\frac{n+4}{n-4}},\, u_\varepsilon
>0\,\,  \mbox{  in  }  A_\varepsilon; \,
u_\varepsilon= \Delta u_\varepsilon= 0\,\, \mbox{ on }  \partial
A_\varepsilon $, where $\{A_\e \subset \R^n,\, \e >0 \}$ is a
family of bounded annulus shaped domains such that $A_\e$ becomes
``thin'' as $\e\to 0$. Our main result is the following: Assume
$n\geq 6$ and let $C>0$ be a constant. Then  there exists
$\varepsilon_0>0$ such that for any $\varepsilon <\varepsilon_0$,
the problem
$ P_{\varepsilon}$ has no solution $u_\e$, whose energy,
$\int_{A_\varepsilon}|\Delta u_\varepsilon |^2$ is less than $C$.
Our proof involves rather delicate analysis of asymptotic profiles
of solutions $u_\e$ when $\e\to 0$
}
\section{Introduction }
Let us consider the following nonlinear elliptic problem under the Navier boundary condition
 $$
P(\Omega )\quad \left\{
\begin{array}{ccccc}
 \Delta^2 u &=& u^{q}, & u >0 & \mbox{  in  }\,  \Omega\\
 u &=&\Delta u=0 & & \mbox{ on }\,  \partial  \Omega,
\end{array}
\right.
$$
\noindent
where $\Omega$ is a bounded regular domain in $\R^n,\,  n\geq 5$ and $q+1 =
2n/(n-4)$ is the critical Sobolev exponent for  the embedding $H^2(\O)\cap H_0^1(\O)$ into $L^{q+1}(\O)$.\\
The interest in this type of equation comes from the fact that it resembles  some
geometric equations involving Paneitz operator (see for instance \cite{C} and \cite{DHL}).\\
It is  well known that if $\Omega$ is starshaped, $P(\Omega )$ has no
solution (see Mitidieri \cite{M} and Van der Vost \cite{V}) and if $\Omega$ has nontrivial topology,
in the sense that  $H_{k}(\Omega ;Z/2Z) \neq
0$ for some $k\in\N,$ Ebobisse and Ould Ahmedou \cite{EO} have shown that  $P(\Omega )
\,$ has a solution. Nevertheless,  Gazzola, Grunan and Squassina \cite{GGS}
gave
the example of contractible domain on which  a solution still exists, showing that both topology and geometry of the domain play a role.\\
In contrast with the subcritical case $q<\frac{n+4}{n-4}$, the associated variational problem happens to be lacking of compactness, that is the functional corresponding to $P(\O)$ does not satisfy the Palais-Smale condition. This means that there exist sequences along which the functional is bounded, its gradient goes to zero, and which do not converge.  Such a fact follows fron the noncompactness of the embedding of  $H^2(\O)\cap H_0^1(\O)$ into $L^{q+1}(\O)$. Since this lack of compactness, the standard variational techniques do not apply and therefore
the question related to existence or nonexistence of solutions of $P(\Omega )
$ remained open.\\
In this paper, we study the problem $P(\Omega)$  when $\Omega  = A_{
\varepsilon}$ is a ringshaped open set in $\R^n$and $\varepsilon \rightarrow
0$. More precisely, let $f$ be any smooth function :
$$
f:\R^{n-1}\longrightarrow [1,2]\, \, ,
(\theta_1 ,...,\theta_{n-1} ) \longrightarrow f(\theta_1 ,...,\theta_{n-1})
$$
\noindent
which is periodic of period $\pi$ with respect to  $\theta_1 ,. . . ,\theta_{
n-2}$ and of period $2\pi$ with respect to $\theta_{n-1}$.\\
We set
$$
S_1 (f) = \left\{ x\in\R^n / r= f(\theta_1 ,. . . ,\theta_{n-1})  \right\},
$$
\noindent
where $ (r ,\theta_1 ,. . . ,\theta_{n-1})$ are the polar coordinates  of  $x$.\\
For $ \e $ positive small enough, we introduce the following map
$$
g_\e:S_1 (f)\longrightarrow g_\e (S_1 (f)) = S_2(f), \, \, \, x \longmapsto g_\e (x)= x + \e n_x,
$$
where $ n_x $ is the outward normal to $ S_1 (f) $ at $ x $.
We denote by $ (A_\e)_{\e > 0} $ the family of annulus shaped domains in $ \R ^n $ such that $ \partial A_\e = S_1 (f) \cup S_2 (f) $.\\
Our main result is the following Theorem.
\begin{thm}\label{t:11}
Assume that $n\geq 6$ and let C be any positive constant. Then there exists $\varepsilon_0 >0$  such that, for any $\varepsilon <\varepsilon_0$,  the problem
$
P_{\varepsilon}:\,
\Delta^2 u_\varepsilon =  u_\varepsilon^{\frac{n+4}{n-4}},\, u_\varepsilon >0
 \mbox{  in  }  A_\varepsilon,\,\,
u_\varepsilon = \D u_\e =0 \mbox{ on }  \partial A_\varepsilon
$,
has no solution such that $\int_{A_\varepsilon}|\D u_\varepsilon |^2 \leq C.$
\end{thm}
\begin{rem}
 We believe the result to be true also for n = 5 (see Remark \ref{r:12} below).
\end{rem}
The proof of Theorem \ref{t:11} is inspired by strong arguments which we developed for the corresponding second order equation \cite{BEH}. It involves rather delicate analysis of asymptotic profiles of solutions when $\e$ tends to zero. Compared with the second order case, further technical problems have to be solved by means of delicate and careful estimates.\\
The plan of the present paper is as follows. In section 2, arguing by contradiction, we suppose that $(P_\e)$ has a solution $u_\e$ with a bounded energy and we study the asymptotic behavior of such a solution, we prove that $u_\e$ blows up at finite points. Then we give in section 3 the characterization of blow up points. Lastly, section 4 is devoted to the proof of our theorem.

\section{Asymptotic  Behavior of Bounded Energy Solutions}
In this section we  suppose that $P_{\varepsilon}$ has a solution $u_{\varepsilon}$ which satisfies $\int_{A_{\varepsilon}}|\D u_{\varepsilon}|^2\leq C ,$ C being a given constant. Our purpose is to study  the asymptotic behavior of $u_{\varepsilon}$ when $\varepsilon $ tends to zero. We prove that $u_{\varepsilon}$ blows up at  p points ($p\in\N^*$).\\
In order to formulate the result of this section, we need to fix some notation.\\
We denote by $G_{\varepsilon}$ the Green's function of $\D^2$
defined by: $\forall\, x\in A_{\varepsilon }$
\begin{eqnarray}\label{e:11}
\Delta^2 G_{\varepsilon}(x,.) = c_n\delta_x\, \mbox{  in  } \,A_{\varepsilon},\quad
\quad G_{\varepsilon}(x,.) = \D G_\e(x,.)=0 \,\mbox{ on } \, \partial  A_{\varepsilon},
\end{eqnarray}
 \noindent
where $\delta_x$ is the Dirac mass at $x$ and $c_n = (n-4)(n-2)|S^{n-1}|$.\\
We denote by  $H_{\varepsilon}$ the regular part of $G_{\varepsilon},$  that
is,
\begin{eqnarray}\label{e:12}
\quad H_{\varepsilon}(x_1 ,x_2) = {|x_1-x_2|^{4-n}} - G_{\varepsilon}
(x_1 , x_2) ,\mbox{ for } (x_1 ,x_2)\in A_{\varepsilon}\times A_{\varepsilon}.
\end{eqnarray}
 \noindent
For $p\in\N^*$ and $ x= (x_1 ,...,x_p)\in A_{\varepsilon}^{p} $,
we denote by $M = M_{\varepsilon}(x)$ the matrix defined by
\begin{eqnarray}\label{e:13}
  M = (m_{ij})_{1\leq i,j\leq p} ,\mbox{ where  } m_{ii} = H_{\varepsilon}(x_i ,x_i) , m_{ij} = -G_{\varepsilon}(x_i ,x_j) ,i\neq j
\end{eqnarray}
\noindent
and define $\rho_{\varepsilon}(x) $ as the least eigenvalue of $M(x)$  ($ \rho_{\varepsilon}(x) = -\infty $ if $x_i = x_j$ for some  $i\neq j$). \\
For $a\in\R^n$ and $\lambda >0,$  $\delta_{(a,\lambda )}$ denotes the function
\begin{eqnarray}\label{e:14}
 \delta _ {(a,\lambda )}(x) = c_0 \lambda
 ^{\frac{n-4}{2}}(1+\lambda^2|x-a|^2)^{\frac{4-n}{2}}.
\end{eqnarray}
\noindent
It is well known (see \cite{Lin}) that if $c_0$ is suitably chosen
$ (c_0 = ((n-4)(n-2)n(n+2))^{\frac{n-4}{8}}) $ the function
$\delta _{(a,\lambda )}$ are the only solutions of equation
\begin{eqnarray}\label{e:15}
\quad \Delta^2 u =  u^{\frac{n+4}{n-4}} , u>0 \mbox{  in  } \R^n
\end{eqnarray}
\noindent
and they are also the only minimizers for the Sobolev inequality, that is
\begin{eqnarray}\label{e:16}
 S =\inf\{|\D u|^{2}_{L^2(\R^n)}|u|^{-2}_{L^{\frac{2n}{n-4}}(\R^n)}
,\, s.t.\, \D u\in L^2 ,u\in L^{\frac{2n}{n-4}} ,u\neq 0 \}.
\end{eqnarray}
\noindent
We also denote by $P_{\varepsilon}\delta_{(a,\lambda)}$
the projection of $\delta_{(a,\lambda )}$ on
$H^2(A_\e)\cap H_0^1(A_{\varepsilon}),$  that is,
$$
\Delta^2 P_{\varepsilon}\delta_{(a,\lambda )} = \Delta^2
\delta_{(a,\lambda )} \mbox{  in  } A_{\varepsilon}, \,\D
P_\e\d_{(a,\l)}= P_{\varepsilon}\delta_{(a,\lambda )}  = 0 \mbox{
on } \partial A_{\varepsilon}.
$$
Lastly, let $E_\e:=  H^2(A_\e)\cap H_0^1(A_{\varepsilon})$
equipped with the norm $||.||$ and the corresponding inner product
$(.,.)$ defined by
\begin{eqnarray}\label{e:19}
||u||&=\left(\int_{A_\e}|\D u|^2\right)^{1/2},\quad
(u,v)&=\int_{A_\e} \D u\D v,\qquad u,v \in E_\e
\end{eqnarray}
  and we define on $E_\e\setminus \left\{0\right\}$ the functional
\begin{eqnarray}\label{e:17}
 J_{\varepsilon}(u) = \left(\int_{A_{\varepsilon}}|\D u|^2\right)\left(
 \int_{A_{\varepsilon}}|u|^\frac{2n}{n-4} \right)^{\frac{4-n}{n}}
\end{eqnarray}
\noindent
whose positive critical points, up a multiplicative constant, are
solutions  of $P_\varepsilon$.\\
Now  we are able to state the main result of this section.
\begin{thm} \label{t:21}
Let $u_{\varepsilon}$ be a solution of problem $P_\varepsilon,$  assume
$
\int_{A_{\varepsilon}}|\D u_\e|^2 \leq C,
$
where $C$ is a positive constant independent of $\varepsilon.$
Then, after passing to a subsequence, there exist  $p\in\N^*$, $(a_{1,
\varepsilon},...,a_{p,\varepsilon})\in A_\varepsilon^p$  and $(\lambda_{1,
\varepsilon},...,\lambda_{p, \varepsilon})\in (\R_+^*)^p$ such that
$$
\bigg|\bigg|{ u_\varepsilon - \sum_{i=1}^{p} P_{\varepsilon}
\delta_{(a_{i , \varepsilon},\lambda_{i ,\varepsilon})}
}\bigg|\bigg| \rightarrow 0, \lambda_{i,\varepsilon}\rightarrow
+\infty , \lambda_{i,\varepsilon}d_{i,\varepsilon}\rightarrow
+\infty, \varepsilon_{ij} \rightarrow 0 \,\,\mbox{ as }\,\,
\e\rightarrow 0,
$$
where $ d_{i,\varepsilon}=d(a_{i,\varepsilon},\partial
A_\varepsilon )$ and where $ \varepsilon_{ij} = (
\frac{\lambda_{i,\varepsilon}}{\lambda_{j,\varepsilon}} +
\frac{\lambda_{j,\varepsilon}}{\lambda_{i,\varepsilon}}
+\lambda_{i, \varepsilon}\lambda_{j,\varepsilon}|a_{i,\varepsilon}
- a_{j,\varepsilon} |^2)^{-\frac{n-4}{2}} $.
\end{thm}

The goal of this section is to prove  Theorem \ref{t:21}. To this aim, we
begin by proving the following lemma :
\begin{lem}\label{l:21}
We have the following claims
\begin{align*}
1.\quad &\int_{A_\varepsilon}|\D u_\varepsilon |^2 \not\longrightarrow 0,\, \mbox {when} \, \,  \varepsilon \longrightarrow 0.\\
2.\quad &
M_\varepsilon \longrightarrow + \infty,\, \mbox {when} \, \,  \varepsilon \longrightarrow 0, \, \mbox{where}\, \,  M_\e =|u_{\varepsilon}|_{L^{\infty}(A_\varepsilon)}.\\
3.\quad & \exists \,c>0 \mbox{ such that for }\, \e\, \mbox{ small enough, we have }\, \e M_\e^{\frac{2}{n-4}} \geq c.
\end{align*}
\end{lem}
\begin{pf}
Since $u_\varepsilon$ is a solution of $P_\varepsilon$, it is clear that
$$
\left(  \int_{A_\varepsilon}| u_\varepsilon|^{\frac{2n}{n-4}} \right)^{\frac{n-4}{n}}\leq\frac{1}{S}\int_{A_\varepsilon}|\D u_\varepsilon |^2 = \frac{1}{S}\int_{A_\varepsilon}u_\varepsilon^{\frac{2n}{n-4}},
$$
where S denotes the Sobolev constant defined in \eqref{e:16}. Thus
 $$
S_n := S^\frac{n}{4} \leq \int_{A_\varepsilon}|\D u_\varepsilon
|^2=\int_{A_\varepsilon}u_{\varepsilon}^{\frac {2n}{n-4}}\leq c \,
 \varepsilon \, M_\varepsilon^{\frac {2n}{n-4}}.
$$
Therefore Claims 1 and 2 are proved. To prove Claim 3, we observe that
$$
\int_{A_\varepsilon}| \D {u_\varepsilon}| ^2 = \int_{
A_\varepsilon}u_{\varepsilon}^\frac {2n}{n-4}
\leq M_\varepsilon^{\frac {8}{n-4}}\int_{
A_\varepsilon}u_\varepsilon^2(x)dx.
$$
Now it is clear that
$$
\int_{A_\varepsilon}u_\varepsilon^2(x)dx =\varepsilon^n \int_{
B_\varepsilon}v_\varepsilon^2(X)dX,
$$
where $ v_\varepsilon(X)=u_\varepsilon (\varepsilon X)$
and where $  B_\varepsilon=\varphi(A_\varepsilon),$
with $  \varphi: x \longmapsto \varphi (x)={\varepsilon}^{-1}x.$\\
Notice that
$$
\varepsilon^n \int_{B_\varepsilon}v_\varepsilon^2(X)dX \leq \frac {\varepsilon^n}{c_\varepsilon}\int_{B_\varepsilon}|
 \D v_\varepsilon (X)|^2 dX = \frac{\varepsilon^4}
{c_\varepsilon} \int_{A_\varepsilon}|\D u_\varepsilon(x)|^2dx,
$$
where $c_\e$ is the least eigenvalue of $\D^2$ on $B_\e$ with the Navier boundary condition.\\
Thus
$$
c_\varepsilon \leq \varepsilon^4 M_\varepsilon^{8/(n-4)}.
$$
As in Lin \cite{L}, we can  prove that $\lim_{\varepsilon \to 0}
c_\varepsilon = c >0$, therefore Claim 3 holds and our lemma follows.
\end{pf} \\
\noindent
Now let $\widetilde A_\varepsilon =  M_\varepsilon^\frac
{2}{n-4}(A_\varepsilon - a_{1,\varepsilon }),$ where $a_{1,
\varepsilon}\in A_\varepsilon$ such that $ M_\varepsilon =
 u_\varepsilon(a_{1,\varepsilon}),$ and we denote by
$v_\varepsilon$ the rescaled function defined on $\widetilde A_\varepsilon$ by
\begin{eqnarray}\label{e:21}
v_\varepsilon(X) = M_{\varepsilon}^{-1}u_\varepsilon (a_{1,
\varepsilon} + M_{\varepsilon}^\frac{-2}{n-4} X).
\end{eqnarray}
\noindent
It is easy to see that $v_\varepsilon$ satisfies
\begin{eqnarray}\label{e:22}
   \left\{
\begin{array}{ccccc}
\Delta^2 v_\varepsilon&=& v_\varepsilon^{\frac{n+4}{n-4}},&
0< v_\varepsilon\leq 1 &\mbox{  in   } \widetilde{ A}_\varepsilon \\
v_\varepsilon(0)&=&1,&  v_\varepsilon =\D v_\e=0 &\mbox{ on }
 \partial\widetilde{  A}_\varepsilon.
\end{array}
\right .
\end{eqnarray}
\noindent
Observe that
$$
\int_{\widetilde{A}_\varepsilon}|\D v_\varepsilon |^2 =
\int_{A_\varepsilon}|\D u_\varepsilon |^2 =\int_{\widetilde{
A}_\varepsilon}v_\varepsilon^{\frac{2n}{n-4}}  =\int_{A_\varepsilon}
u_\varepsilon^{\frac{2n}{n-4}} \leq C.
$$
\noindent
Let us distinguish the following cases:\\
1. \, $M_\varepsilon^\frac{2}{n-4}d(a_{1,\varepsilon} ,
\partial A_\varepsilon )\rightarrow +\infty , \, \mbox{ when }\, \varepsilon
\rightarrow 0.$\\
2. \, $M_\varepsilon^\frac{2}{n-4}d(a_{1,\varepsilon} ,
\partial A_\varepsilon )$ tends to 0.\\
3. \, $M_\varepsilon^\frac{2}{n-4}d(a_{1,\varepsilon} ,
\partial A_\varepsilon )$ is bounded below and above.\\
As in the proof of Lemma 2.3 \cite{BEH}, we can show that case 2 cannot occur. Now we are going to prove that case 3 cannot also occur. Arguing by contradiction, let us suppose that case 3 occurs. Then it follows from \eqref {e:22} and
standard elliptic theories that there exists some positive
function $v$, such that ( after passing to subsequence ), $ v_\e
\to v $ in $ C^2_{loc} ( \O ) $, where $ \O $ is a half space or a
strip of $ \R^n $, and $v$ satisfies
 \begin{eqnarray}\label{o}
   \left\{
\begin{array}{cccccc}
\Delta^2 v& =& v^{\frac{n+4}{n-4}} , &0\leq v \leq 1  &\mbox{  in  }&  \Omega \\
v(0)& =& 1,&  v =\D v=0& \mbox{ on }&  \partial  \Omega.
\end{array}
\right .
\end{eqnarray}
\noindent But if $\Omega $ is a half space or a strip of $\R^n$,
then $v$ must vanish identically (see \cite{M}). Thus we derive a
contradiction. So we are in the first
case and therefore there exists some positive function $v,$
such that (after passing to a subsequence), $v_\varepsilon
\longrightarrow v$  in $C_{loc}^2(\R^n),$ and $v$ satisfies
\eqref{o} with $\O=\R^n$ and without boundary conditions. It
follows from Lin \cite{Lin}
$$
 v(X) =  \delta_{(0,\a _n )}(X), \quad \mbox { with } \a _n
 = ((n-4)(n-2)n(n+2))^{ {-1}/{4}}.
$$
Hence
$$
M_\varepsilon ^{-1}u_{\varepsilon}(a_{1,\varepsilon} +
M_\varepsilon ^{\frac{-2}{n-4}} X ) -  \delta_{(0,\a _n )}
(X)\rightarrow 0\mbox{ in  }C_{loc}^2(\R^n) ,\mbox{  when
}\varepsilon \rightarrow 0.
$$
\noindent
Observe that
$$
M_\varepsilon ^{-1} u_{\varepsilon}(a_{1,\varepsilon} +
M_\varepsilon ^{\frac{-2}{n-4}} X ) -  \delta_{(0,\a _n )}(X) =
M_\varepsilon ^{-1}\left(u_\varepsilon(x) -  \delta_{(a_{1,
\varepsilon}, \lambda_{1,\varepsilon})}(x)\right),\,\,\mbox{
with }\,\,\lambda_{1,\varepsilon}=\a _n M_\varepsilon ^{\frac{2}{n-4}}.
$$
In the sequel, we denote by $u_\varepsilon^1$ the function defined on
$A_\varepsilon$ by
\begin{eqnarray}\label{e:24}
  u_\varepsilon^1(x) = u_\varepsilon (x) - P_\varepsilon
\delta_{(a_{1,\varepsilon} ,\lambda_{1,\varepsilon})}(x).
\end{eqnarray}
\noindent
Notice that   $\lambda_{1,\varepsilon}\rightarrow  +\infty$ and
$\lambda_{1,\varepsilon}d(a_{1,\varepsilon} ,
\partial A_\varepsilon ) \rightarrow +\infty$ when $\varepsilon\rightarrow 0$.\\
Now we need to prove the following lemma :
\begin{lem}\label{l:24}
Let $u_\varepsilon^1 $ be defined by \eqref{e:24}. Then we have \\
{\bf i /}$ \Delta^2 u_\varepsilon^1 = |u_\varepsilon^1
|^{\frac{8}{n-4}}u_\varepsilon^1 +g_\varepsilon ,$ with
$|g_\varepsilon |_{H^{-2}_{(A_\varepsilon )}}\rightarrow 0$ ,
when $\varepsilon\rightarrow 0$.\\
{\bf ii /} $\int_{A_\varepsilon}|\D u_\varepsilon^1|^2 =
 \int_{A_\varepsilon}|\D u_\varepsilon|^2 - S_n + o(1)$.\\
{\bf iii /}$  \int_{A_\varepsilon}|u_\varepsilon^1|^{\frac{2n}{n-4}} =   \int_{A_\varepsilon}|u_\varepsilon|^{\frac{2n}{n-4}} -
S_n + o(1)$, where $ S_n = S^{\frac{n}{4}}$.
\end{lem}
\begin{pf}
{\bf i/} We have
$$
 \Delta^2 u_\varepsilon^1
=u_\varepsilon^{\frac{n+4}{n-4}} - \delta_{(a_{1 , \varepsilon} ,
\lambda_{1, \varepsilon})}^{\frac{n+4}{n-4}}
=|u_\varepsilon^1|^{\frac{8}{n-4}}u_\varepsilon^1
+g_\varepsilon,\,\,\mbox{with}\,\,g_\varepsilon
=u_\varepsilon^{\frac{n+4}{n-4}} - \delta_{(a_{1 ,\varepsilon} ,
\lambda_{1, \varepsilon})}^{\frac{n+4}{n-4}}
-|u_\varepsilon^1|^{\frac{8}{n-4}}u_{\varepsilon}^{1}.
$$
Observe that
$$
g_\varepsilon = O\left( P_\varepsilon\delta ^{\frac{8}{n-4}} |
u_\varepsilon - P_\varepsilon\delta | +  | u_\varepsilon -
P_\varepsilon\delta
 |^{\frac{8}{n-4}}P_\varepsilon\delta \right) + O\left(
\delta_\varepsilon^{\frac{8}{n-4}}\th_\e\right),
$$
\noindent
where $P_\varepsilon\delta = P_\varepsilon \delta_{(a_{1 ,
\varepsilon} , \lambda_{1, \varepsilon})}$, $\delta_
\varepsilon = \delta_{(a_{1 ,\varepsilon} , \lambda_{1,
\varepsilon})}$ and $\th_\e=\d_\e-P_\e\d$.\\
Since $L^{\frac{2n}{n+4}}\hookrightarrow H^{-2} $,
it is sufficient to prove that
$$
\int_{A_\varepsilon}\delta_\varepsilon ^{\frac{16n}{n^2-16}}
| u_\varepsilon - P_\varepsilon\delta |^{\frac{2n}{n+4}}\rightarrow 0
 \, \mbox{and}  \int_{A_\varepsilon}\delta_\varepsilon ^{\frac{2n}{n+4}}
 | u_\varepsilon - P_\varepsilon\delta |^{\frac{16n}{n^2-16}}\rightarrow 0,
  \mbox{when } \varepsilon \rightarrow  o.
$$
\noindent
Observe that
\begin{align*}
\int_{A_\varepsilon}\delta_\varepsilon ^{\frac{16n}{n^2-16}}&
 | u_\varepsilon - P_\varepsilon\delta |^{\frac{2n}{n+4}}
\leq c \int_{A_\varepsilon}\delta_\varepsilon ^{\frac{16n}{n^2-16}}
| u_\varepsilon - \delta_\varepsilon |^{\frac{2n}{n+4}}  + c
\int_{A_\varepsilon}\delta_\varepsilon ^{\frac{16n}{n^2-16}}
  \th _\varepsilon ^{\frac{2n}{n+4}}\\
& \leq c\int_{\widetilde{A}_\varepsilon} \delta_{(0, \a _n)}
^{\frac{16n}{n^2-16}}(X) | v_\varepsilon  (X) -
M_{\varepsilon}^{-1}\delta_\varepsilon (a_{1,\varepsilon} +
M_{\varepsilon}^{\frac{-2}{n-4}} X) |^{\frac{2n}{n+4}}dX + O\left(
|\th_\varepsilon |^{\frac{2n}{n+4}}_{L^{\frac{2n}{n-4}} } \right).
\end{align*}
The function $\theta_\e$ satisfies (see \cite{BH1})
\begin{eqnarray}\label{theta}
|\theta_\e|_{L^{2n/(n-4)}}=O\left((\l_{1,\e}d_{1,\e})^{(4-n)/2}\right)=o(1).
\end{eqnarray}
Regarding  the first term, let $R$ be a large constant such that
$\int_{\R^n\setminus B(0,R)}\delta _{(0,\a _n)}^{\frac{2n}{n-4}}=
o(1)$. Then, using the Holder's inequality and the fact that
$\int_{\widetilde{A}_\e} v_\e ^{{2n}/(n-4)}\leq C$, we derive that
\begin{eqnarray}
\int_{\widetilde{A}_\varepsilon\setminus B(0,R)}\delta_{(0, \a _n
 )} ^{\frac{16n}{n^2-16}} (X)| v_\varepsilon  (X) - \delta_{(0,\a _n )}( X)
|^{\frac{2n}{n+4}}dX \leq c_1\left( \int_{\R^n\setminus B(0,R)}
\delta_{(0, \a _n )} ^{\frac{2n}{n-4}}\right)^{\frac{8}{n+4}}=o(1)
\end{eqnarray}
Now, we need to estimate the following integral
$$
\int_{ B(0,R)}\delta_{(0,\a _n )} ^{\frac{16n}{n^2-16}}(X)
| v_\varepsilon  (X) - \delta_{(0,\a _n )}( X)
|^{\frac{2n}{n+4}}dX \leq C\int_{ B(0,R)}
 | v_\varepsilon  (X) - \delta_{(0,\a _n )}( X)
|^{\frac{2n}{n+4}}dX\rightarrow   0,
$$  when $\varepsilon\rightarrow 0$,
indeed $v_\varepsilon  - \delta_{(0,\a _n )} \longrightarrow 0$ in $ C_{loc}^2(\R^n).$ \\
In the same way, we prove that
$$
\int_{A_\varepsilon}\delta_\varepsilon ^{\frac{2n}{n+4}}
| u_\varepsilon - P_\varepsilon\delta |^{\frac{16n}{n^2-16}}\longrightarrow 0,
 \mbox{when}\, \,  \varepsilon \longrightarrow 0.
$$
\noindent
{\bf ii/} We also have
$$
\int_{A_\varepsilon}|\D u_\varepsilon^1|^2 = \int_{A_\varepsilon}
|\D u_\varepsilon|^2 + \int_{A_\varepsilon}|
\D P_\varepsilon\delta |^2 -2\int_{A_\varepsilon}
\D u_\varepsilon\D P_\varepsilon\delta.
$$
\noindent
Observe that
\begin{eqnarray*}
\int_{A_\varepsilon}|\D P_\varepsilon\delta|^2& =
&\int_{A_\varepsilon}\delta_\varepsilon^{\frac{n+4}{n-4}}
P_\varepsilon\delta
=\int_{A_\varepsilon}\delta_\varepsilon^{\frac{2n}{n-4}} -
\int_{{A}_\varepsilon}\delta_\varepsilon^{\frac{n+4}{n-4}}\th_\varepsilon \nonumber\\
&=&\int_{\widetilde{A_\varepsilon}}\delta_{(0, \a _n)}^{\frac{2n}{n-4}}-
\int_{A_\varepsilon}\delta_\varepsilon^{\frac{n+4}{n-4}}\th_\varepsilon
=S_n  - \int_{\R^n\setminus\widetilde{A}_\varepsilon}
\delta_{(0,\a _n )}^{\frac{2n}{n-4}} - \int_{A_\varepsilon}
\delta_\varepsilon^{\frac{n+4}{n-4}}\th_\varepsilon.
\end{eqnarray*}
\noindent
For the 2nd integral, we have
$$
\int_{A_\varepsilon}\delta_\varepsilon^{\frac{n+4}{n-4}}
\th_\varepsilon \leq C |
\th_\varepsilon |_{L^{\frac{2n}{n-4}}(A_\varepsilon)}\leq
{c}{(\lambda_{1,\varepsilon}\, d(a_{1,\varepsilon},
\partial{ A_\varepsilon}) )^{\frac{4-n}{2}}}\longrightarrow 0
\mbox{ when }\varepsilon\longrightarrow 0,
$$
\noindent
For the first integral, we have
$$
\int_{\R^n\setminus\widetilde{A}_\varepsilon}
\delta_{(0,\a _n)}^{\frac{2n}{n-4}} =o(1)
$$
\noindent
indeed $\widetilde{A}_\varepsilon\longrightarrow \R^n$ and
 $\delta_{(0,\a _n)}\in L^{\frac{2n}{n-4}}(\R ^n)$. Then
$$
\int_{A_\varepsilon}|\D  P_\varepsilon\delta|^2 = S_n + o(1).
$$
\noindent
We also have
$$
\int_{A_\varepsilon}\D u_\varepsilon \D P_{\varepsilon} \delta
= \int_{A_\varepsilon}\D ( u_\varepsilon  - P_\varepsilon\delta )
\D  P_\varepsilon\delta + \int_{A_\varepsilon}|
\D  P_\varepsilon\delta|^2.
$$
\noindent
Observe that
\begin{align*}
\int_{A_\varepsilon}\D ( u_\varepsilon  -
P_\varepsilon\delta ) & \D  P_\varepsilon\delta
= \int_{A_\varepsilon}(u_\varepsilon -
P_\varepsilon \d )\delta _{\varepsilon} ^{\frac{n+4}{n-4}}
=\int_{A_\varepsilon}(u_\varepsilon -
\delta_\varepsilon )\delta_\varepsilon^{\frac{n+4}{n-4}}
+ \int_{A_\varepsilon}\th_\varepsilon\delta_\varepsilon ^{\frac{n+4}{n-4}}\nonumber\\
&=\int_{\widetilde{A_\varepsilon}}(v_\varepsilon -
\delta_{(0, \a _n)}) \delta_{(0, \a _n)}^{\frac{n+4}{n-4}} +o(1)\nonumber\\
&\leq  \int_{ B(0,R)}(v_\varepsilon -
\delta_{(0, \a _n)}) \delta_{(0, \a _n)}^{\frac{n+4}{n-4}}
+ \int_{ \R^n\setminus B(0,R)}(v_\varepsilon -
\delta_{(0, \a _n)}) \delta_{(0, \a _n)}^{\frac{n+4}{n-4}} +o(1)
\end{align*}
\noindent Notice that, on one hand
$$
 \int_{ \R^n\setminus B(0,R)}(v_\varepsilon -
\delta_{(0,\a _n)}) \delta_{(0, \a _n)}^{\frac{n+4}{n-4}}
\leq C\left( \int_{ \R^n\setminus B(0,R)}
 \delta_{(0, \a _n)}^{\frac{2n}{n-4}}\right)^{\frac{n+4}{2n}} = o(1).
$$
\noindent
On the other hand
$$
 \int_{ B(0,R)}(v_\varepsilon - \delta_{(0, \a _n)})
\delta_{(0, \a _n)}^{\frac{n+4}{n-4}} = o(1), \mbox{ since }
v_\varepsilon\longrightarrow \delta_{(0, \a _n)} \mbox{ in } C_{loc}^2(\R^n).
$$
\noindent
Then
$$
\int_{A_\varepsilon}\D u_\varepsilon\D  P_\varepsilon\delta =
S_n + o(1).
$$
\noindent
Thus Claim ii/ of Lemma \ref{l:24} follows.\\
{\bf iii/ } The proof of iii/  in Lemma \ref{l:24} is similar to that
of ii/, so we will omit it.
\end{pf} \\
Now we distinguish two cases :\\
{\bf i /} $ \int_{A_\varepsilon}|\D  u_\varepsilon^1
|^2\longrightarrow 0 $ when $\varepsilon\longrightarrow 0$. \\
{\bf ii /} $ \int_{A_\varepsilon}|\D  u_\varepsilon^1
|^2\not\longrightarrow 0 $ when $\varepsilon\longrightarrow 0$. \\
If $ \int_{A_\varepsilon}|\D u_\varepsilon^1|^2\longrightarrow 0$ , the
proof of Theorem \ref{t:21} is finished.\\
In the sequel, we  consider the second case, that is $ \int_{A_\varepsilon}|
\D u_\varepsilon^1|^2\not\longrightarrow 0,  $ when $\varepsilon
\longrightarrow 0$ and we are going to look for a second point of blow up of
 $u_\e$. \\
In order to simplify our notation, in remainder we often omit the index
$ \e $ of $ a_\e $ and $ \l _\e $.\\
Let us introduce the following notation :
\begin{align}
&u_\e (a_2)  := \l _2^{\frac{n-4}{2}} = \max _{(A_\e \diagdown B(a_1, \e ))}
u_\e (x)\label{e:25}\\
&h_\e  = \max _ {B(a_1,2\e)}|x-a_1|^{\frac{n-4}{2}}u_\e (x) = |a_1-a_3|^{
\frac{n-4}{2}}u_\e (a_3) = |a_1-a_3|^{\frac{n-4}{2}}\l _3^{\frac{n-4}{2}}.
\label{e:26}
\end{align}
We distinguish two cases :\\
{\bf {Case 1. }} $ h_\e \to + \infty $ when $ \e \to 0 $.\\
{\bf { Case 2.}} $ h_\e \leq c $, when $ \e \to 0 $.\\
Now we study the first case, that is, $ h_\e \to \infty $ when $ \e \to 0 $.\\
Let
$$
\l _4 = \max (\l _2, \l _3 ) := u_\e ^{\frac{2}{n-4}}(a_4).
$$
For $ X \in B(0, \frac{\l _4}{2}|a_1-a_4|)\cap D_\e $, we set
$$
w_\e(X) = \l _4^{\frac{4-n}{2}}u_\e(a_4 + \l _4^{-1} X ), \quad
\mbox {with }\, D_\e = \l _4 (A_\e - a_4).
$$
It is easy to check the following claims
$$
\l _4 |a_1-a_4| \geq (1/2){\l _3}|a_1-a_3|, \quad \mbox {and }\,  \l _4 \e
\geq (1/2) \l _3 {|a_1-a_3|}.
$$
Thus
$$
\l _4 |a_1-a_4|\to +\infty \quad \mbox {and }\, \l _4 \e \to +\infty \quad
\mbox {as } \e \to 0.
$$
We also have
$$
w_\e (X) \leq c, \quad \foral X \in B(0, (1/2){\l _4}|a_1-a_4|) \cap D_\e.
$$
By an argument similar to the one used after the proof of Lemma
\ref{l:21}, we have\\ $ \l _4 d(a_4 , \partial A_\e)\to +\infty$,
as $\e \to 0$. Thus, there exist $ b \in \R^n $ and $ \l > 0 $
such that\\ $ w_\e \to \d _{(b, \l)} $ in $ C^2_{loc}(\R^n) $.
Therefore we have found a second point of blow up $ \bar {a}_2 $
of $ u_\e $ with the concentration $ \bar {\l }_2 $ in this case
($ \bar {a}_2 = a_4 +  {b}/{\l _4} $ and $ \bar {\l }_2 = \l \l _4 $ ).\\
Next we study the second case, that is, $ h_\e $ remains bounded when $ \e \to
0 $, where $ h_\e $ is defined in \eqref{e:26}. In this case we consider two
subcases.\\
{\bf {i.}}\hskip 0.3cm $ \int _ {B(a_1, 2\e)}|u_\e^1|^{\frac{2n}{n-4}} \to 0 $
as $ \e \to 0 $.\\
{\bf {ii.}}\hskip 0.3cm $ \int _ {B(a_1, 2\e)}|u_\e^1|^{\frac{2n}{n-4}}
\nrightarrow 0 $ as $ \e \to 0 $.\\
Let us consider the first subcase. Clearly, we have
$\int_{A_\e\backslash B(a_1, 2\e) } |u_\e
^1|^{2n/(n-4)}\nrightarrow 0$ and \\ $\int_{A_\e\backslash B(a_1,
2\e) } \d_{(a_1,\l_1)}^{2n/(n-4)}\to 0$. Then, there exists $ c
> 0 $ such that
$$
0 < c \leq  \int _ {A_\e \diagdown B(a_1, 2\e)} u_\e
 ^{\frac{2n}{n-4}}\leq c\l _2^2 \int _{A_\e} u_\e ^2 \leq c(\e \l
_2)^2.
$$
Hence, there exists $ \bar c > 0 $ such that
$$
\l _2 |a_1-a_2| \geq \l _2 \e \geq 2 \bar c.
$$
Now, for $ X \in E_\e = \l _2 ( A_\e - a_2) $, we introduce the following
function
$$
U_\e(X) = \l _2^{\frac{4-n}{2}}u_\e (a_2 + \l _2^{-1}X ).
$$
We also have $ \l _2 d(a_2, \partial A_\e) \to
+\infty $. It is easy to see that $ U_\e $ satisfies
$$
U_\e \leq 1, \quad \mbox {in }\, B(0, (1/2){\l _2}|a_1-a_2|).
$$
Thus there exists $ b \in \R^n $ and $ \l > 0 $ such that $ U_\e \to \d _{(b,
\l )} $ in $ C^2_{loc}(\R^n ) $. Therefore we have also found a second point
of blow up $ \bar {a}_2 $ of $ u_\e $ with the concentration $ \bar {\l }_2 $
 in this case  ($ \bar {a}_2 = a_2 + {b}/{\l _2} $ and $ \bar {\l }_2 = \l
 \l _2 $ ). \\
Now we study the second subcase. To this aim,
we introduce the following function defined on $ F_\e = \e ^{-1}(A_\e - a_1) $
 by
$$
W_\e (X) = \e ^{\frac{n-4}{2}}u_\e ^1 (a_1 + \e X ).
$$
Observe that $ F_\e $ ``converges'' to a strip of $ \R^n $ when $ \e \to 0 $.
We notice that $ W_\e $ satisfies
\begin{eqnarray*}
\left\{
\begin{array}{ccccc}
\D^2 W_\e & = & |W_\e|^{\frac{8}{n-4}}W_\e + f_\e & \mbox {in }& F_\e \\
W_\e & = & \D W_\e=0        & \mbox {on}&\partial F_\e,
\end{array}
\right.
\end{eqnarray*}
with $ |f_\e|_{H^{-2}(F_\e)} \to 0 $ as $ \e \to 0 $. \\
We also have
$$
\int _{B(0,2)\cap F_\e}|W_\e|^{\frac{2n}{n-4}} = \int _{B(a_1,2\e)\cap A_\e}
|u^1_\e|^{\frac{2n}{n-4}} \nrightarrow 0,\, \mbox{ as }\, \e \to 0
\mbox{ and }
\int _{F_\e}|\D W_\e |^2 = \int _{A_\e}|\D u^1_\e|^2 \leq C
$$
It is easy to check that there exists some fixed domain $ F \subset B(0,2)\cap
 F_\e $ such that $ |W_\e |^{\frac{2n}{n-4}} \to 0 $ almost everywhere and
$ |W_\e |^{\frac{2n}{n-4}} \nrightarrow 0 $ in $ L^1(F) $. From
Dunford-Pettis'Lemma (\cite {Be}), we have
\begin{eqnarray}\label{e:27}
\exists \d_0 > 0, \quad \exists \a _\e > 0, \, \a _\e \to 0, \quad \exists
b_\e \in F \quad \mbox {s.t.} \int _{B(b_\e,\a _\e )\cap F_\e }|W_\e |^{
\frac{2n}{n-4}} \geq \d _0.
\end{eqnarray}
We can choose $ b_\e $ and $ \a _\e $ such that $ \a _\e $ is minimum and
$  \int _{B(b_\e,\a _\e )\cap F_\e }|W_\e |^{\frac{2n}{n-4}} = \d _0 $.
\begin{lem}\label{l:25}
Let $ (\a _\e, b_\e )$ be defined by \eqref{e:27} and let $ \bar {\l }_2= (\e
\a _\e )^{-1} $, and $ \bar {a}_2= a_1 + \e b_\e $.
 Then we have
$$
\frac{\lambda_1}{\bar{\lambda}_2}
 \longrightarrow +\infty  \mbox{ or }  \frac{\bar {\lambda}_2}
{\lambda_1} \longrightarrow +\infty  \mbox{ or }
 \lambda_{1} \bar{\lambda}_2|a_{1}
- \bar {a_2} |^2\longrightarrow +\infty
 \mbox{ when } \varepsilon\longrightarrow 0, \,\mbox{ with }\,\lambda_1 = M_\e^{\frac {2}{n-4}}.
$$
\end{lem}
\begin{pf}
We argue by contradiction. Let us suppose that
$ {\lambda _1}/{\bar {\lambda }_2}, $   $ {\bar {\lambda}_2}/{\lambda_1} $
and  $\lambda_1 \bar {\lambda}_2|a_1 - \bar {a}_2 |^2 $
 are bounded when $ \varepsilon \longrightarrow 0 $.\\
For  $X\in\widetilde{A}_\e: = \lambda_1(A_\varepsilon
- a_1 ),$  we introduce $ \o _\varepsilon $ defined by
\begin{eqnarray}\label{e:28}
 \o _\varepsilon(X) = M_\varepsilon ^{-1}u_\varepsilon^1(a_
1 + \lambda_1^{-1}X).
\end{eqnarray}
\noindent Observe that, on one hand
\begin{eqnarray*}
\int_{B(\lambda _1(\bar {a}_2 - a_1), {\lambda _1}/{
\bar {\lambda}_2} )\cap \widetilde{A}_\varepsilon}
|\o _\varepsilon (X)|^{\frac{2n}{n-4}}dX &=& \int_{B(\bar {a}_2 ,
 {1}/{ \bar {\lambda }_2} )\cap {A}_\varepsilon}
|u_\varepsilon^1 (x)|^{\frac{2n}{n-4}}dx\\
& =& \int _{B(b_\e ,\a _\e )\cap F_\e}
|W_\e(X)|^{\frac{2n}{n-4}}dX = \d _0 > 0.
\end{eqnarray*}
\noindent
On the other hand, since $ \lambda _1
|\bar {a}_2- a_1| $\, and
$ {\lambda _1}/{ \bar {\lambda }_2} $ are bounded, we have
$$
\exists R >0 \,  \mbox{such that } \,  \,    B(\lambda _1
(\bar {a}_2 - a_1),  {\lambda _1}/{\bar {\lambda }_2})  \subset B(0, R).
$$
\noindent
Thus
\begin{align*}
 \int_{B(\lambda _1(\bar {a}_2
- a_1), {\lambda _1}/{\bar {\lambda }_2} )\cap
\widetilde{A}_\varepsilon}&|\o _\varepsilon (X)
|^{\frac{2n}{n-4}}dX
\leq \int_{B(0 ,R )\cap{\widetilde{A}_\varepsilon}}
|\o _\varepsilon (X)|^{\frac{2n}{n-4}}dX\nonumber\\
& = \int_{B(0 ,R )\cap{\widetilde{A}_\varepsilon}}
|M_\varepsilon ^{-1}u_\varepsilon (a_1 + \lambda _1^{-1} X) -
M_\varepsilon ^{-1}P_\varepsilon \delta(a_1 +
\lambda _1^{-1} X)|^{\frac{2n}{n-4}}dX\nonumber\\
& \leq  c\int_{B(0,R)}|v_\varepsilon - \delta_{(0,
\a _n)}(X)|^{\frac{2n}{n-4}} + c\int_{A_\varepsilon}|
\th_\varepsilon(x)|^{\frac{2n}{n-4}}.
\end{align*}
\noindent
Therefore
$$
\int_{B(\lambda _1(\bar {a}_2 -
a_1), {\lambda _1}
{\bar {\lambda }_2 } )\cap \widetilde{A}_\varepsilon}
|W_\varepsilon (X)|^{\frac{2n}{n-4}}dX  \longrightarrow 0 \,
\mbox{ when } \varepsilon \longrightarrow 0
$$
\noindent
which yields a contradiction and our Lemma follows .
\end{pf}\\
\noindent
Now we set
$\widetilde{\widetilde{A}}_\varepsilon = \bar {\lambda }_2
(A_\varepsilon - \bar {a}_ 2 )$ and we introduce the function
$V_\varepsilon$ defined by
$$
V_\varepsilon (X) = \bar {\l }_2^{\frac{(4-n)}{2}}u_\varepsilon^1
(\bar {a}_2 +\bar {\lambda}_2^{-1} X ).
$$
Observe that
\begin{eqnarray}\label{e:n}
 \int_{B(0,1)\cap \tilde{\tilde{A}}_\e }
|V_\varepsilon|^{\frac{2n}{n-4}} =  \int_{B(\bar {a}_2 ,
 {1}/{ \bar {\lambda }_2} )\cap {A}_\varepsilon}
|u_\varepsilon^1 |^{\frac{2n}{n-4}} = \d _0 > 0,\, \int |\D V_\e
|^2 \leq C,\, \int |V_\e |^{\frac {2n}{n-4}}\leq C.
\end{eqnarray}
It is easy to see  that there exists some functions V such that
(after passing to a subsequence), $V_\varepsilon \longrightarrow
V$ in $ H_{loc}^2( \O)$ and  V satisfies
\begin{eqnarray} \label{e:29}
\left\{
\begin{array}{cc}
\Delta^2 V  = |V|^{\frac{8}{n-4}}V &\mbox{  in  } \O, \quad V=\D V
=0 \, \, \mbox{on }
  \partial \O \\
\int _{\O}|\D V|^2 \leq C,& \int _{\O}|V|^{\frac {2n}{n-4}} \leq C,
\end{array}
\right .
\end{eqnarray}
where $\O $ is a half space or a strip or a $\R^n$.\\
From \eqref{e:n}, it is easy to see that $ V \ne 0 $.
\begin{lem}\label{l:26}
Let V be defined by \eqref{e:29}. Then  we have  $V\geq 0.$
\end{lem}
\begin{pf}
We have
\begin{eqnarray}\label{e:210}
V_\varepsilon(X)  =  \bar {\l }_2^{\frac{(4-n)}{2}}
\left(u_\varepsilon(\bar{a}_2 +\bar{\l}_2^{-1}X) - \delta_{( a_1,
\lambda_1)}(\bar {a}_2 +\bar {\lambda }_2^{-1} X )
 + \th_{( a_1,\lambda_1)}(\bar{a}_2
+\bar {\lambda}_2^{-1} X )\right).
\end{eqnarray}
\noindent
Thus, it is sufficient to prove that
$$
 \bar {\l }_2^{\frac{(4-n)}{2}} \delta_{( a_1,
\lambda_1)}(\bar {a}_2 + \bar {\lambda }_2^{-1} X )\longrightarrow
0\mbox{ in }H_{loc}^2(\R^n).
$$
\noindent
Observe that
\begin{eqnarray}\label{e:211}
 I_\varepsilon & := & \int_{B(0,R)}\left( \bar {\l }_2^{\frac{(4-n)}{2}}
\frac{ \lambda_1^{\frac{n-4}{2}}}{(1+
\lambda_1^2|\bar {a}_2  +
\bar {\lambda }_2 ^{-1} X -  a_1
 |^2)^{\frac{n-4}{2}}}
    \right)^{(n+4)/(n-4)}dX \nonumber \\
& = & (\frac{\lambda _1}{\bar {\lambda }_2 })^{(n+4)/2}\int_{B(0,R)}
\frac {dX}
 { \left ( 1 + \left ( {\lambda _1}/
{\bar {\lambda }_ 2} \right)^2 |X- \bar {\lambda }_2
( a_1- \bar {a}_2 )|^2 \right )^{ (n+4)/2}}.
\end{eqnarray}
\noindent
If $ {\lambda _1}/{\bar {\lambda }_2}\longrightarrow 0$ it is clear that
$I_\varepsilon\longrightarrow 0$ when $\varepsilon\longrightarrow 0.$
If ${ \lambda _1}/{\bar {\lambda }_2} \longrightarrow + \infty,$ \, let $y =
 ({\lambda _1}/{\bar {\lambda }_2 })X $. Thus
$$
I_\varepsilon \leq  (\frac{\bar {\lambda }_2}
{\lambda _1})^{(n-4)/2}\int_{\R^n}
\left( \frac{ 1}{1 + |y - \lambda _1
(a_1- \bar {a}_2 )|^2} \right)^{(n+4)/
2} dy \longrightarrow 0 \, \mbox{ as }\,  \varepsilon \to 0.
$$
\noindent
Lastly, if ${\lambda_1}/{\bar {\lambda}_2}
\not\longrightarrow + \infty$ \,   and ${\bar{\lambda}_2}/{\lambda_1}
\not\longrightarrow + \infty$, \,  then, by Lemma \ref{l:25}, we have
$$
\lambda_1 \bar { \lambda }_2
 |a_1- \bar {a}_ 2|^2 \longrightarrow +
\infty, \, \mbox{ when} \, \, \varepsilon \to 0.
$$
\noindent
Observe that for $X\in B(0,R),$ we have
$$
| X- \bar {\lambda }_2 ( a_1 -
\bar { a}_2 )|\geq | \bar {\lambda }_2
( a_ 1- \bar{a}_2)| - |X|\geq c|
\bar {\lambda }_2 ( a_1- \bar {a}_2)|.
$$
\noindent
Therefore
$$
I_\varepsilon \leq {(\frac{\lambda _1}
{\bar {\lambda}_2})}^{\frac{n+4}{2}}\int _{B(0,R)}
{\frac{1}{\left(c(\frac{\lambda _1}{
\bar {\lambda }_2})^2|\bar {\lambda }_2
(a_1 - \bar {a}_2)|^2  \right)^{\frac{n+4}{2}}}}
\leq {\frac{C(R)}{(\lambda _1 \bar {\lambda}_2|a_1- \bar {a}_2|^2)^{
\frac{n+4}{2}}}}\longrightarrow 0, \, \mbox{as }\varepsilon\to 0.
$$
\noindent
Then our Lemma follows.
\end{pf}\\
\noindent
Now, from  \cite{M}, we derive that $\O = \R^n $. Thus,
using \eqref{e:29} and Lemma \ref{l:26}, we also obtain a second point of
blow up of $u_\e $ in this case. Thus in all cases we have built a second
point $a_{2,\e} $ of blow up of $u_\e $ with the concentration $ \l _{2, \e}
$ such that $ \l _{2,\e} \to + \infty $ and $ \l _{2,\e}d(a_{2,\e},\partial
A_\e ) \to
+\infty $ as $ \e \to 0 $. It is clear that  we can proceed by inductions.
Thus we obtain a sequence $(u_\varepsilon^k)_k$  such that
$$
\int_{A_\varepsilon} |\D u_\varepsilon^k|^2 = \int_{A_\varepsilon}
 |\D u_\varepsilon^{k-1}|^2 - S_n +o(1) = \int_{A_\varepsilon}
|\D u_\varepsilon|^2 -kS_n +o(1).
$$
\noindent
Thus
\begin{eqnarray}\label{e:212}
0 \leq \int_{A_\varepsilon} |\D u_\varepsilon^k|^2 =
\int_{A_\varepsilon} |\D u_\varepsilon|^2 -
kS_n +o(1) \leq C- kS_n +o(1).
\end{eqnarray}
\noindent
Since the later term  in \eqref{e:212} will be negative for large  k, the
induction will terminate after some index $p\in \N ^*.$ Moreover, for this
index, we obtain desired claims in  Theorem \ref{t:21}.
\section{ Location of Blow up Points}
\mbox{}
In this section, we give the characterization of blow
up points which we found in section 2. Namely, we prove
the following crucial result :
\begin{thm} \label{t:31}
Let $a_{1,\e}$,..., $a_{p,\e}$ be the points given by
Theorem \ref{t:21}. Then we have $ p\geq 2$. Moreover, if $n\geq 6,$  we have :
$
\exists \, k \leq p , \exists \, i_1 ,...,i_k \in\left\{ 1, 2,...,p \right\}$
  such that
$$
d^{n-4}\rho_\varepsilon (a_{i_1,\varepsilon} ,...,a_{i_k,\varepsilon} )
\rightarrow 0, \quad d^{n-3}{\nabla{\rho}}_{\varepsilon} (a_{i_1,\varepsilon}
,...,a_{i_k,\varepsilon} ) \rightarrow 0.
$$
 In addition, we have $\forall m,l \in \left\{ 1,...,k \right\}$ $|a_{i_m,\varepsilon} -
a_{i_l,\varepsilon} |\leq C_0 d$, where \\ $d =
\min\left\{d(a_{i_l,\varepsilon}) ,\partial A_\varepsilon)\, /
1\leq l\leq k\right\}$ and $C_0$ is a positive constant
independent of $\varepsilon$.
\end{thm}
\begin{rem}\label{r:12}
We believe the result of Theorem \ref{t:31} to be true for
$n=5.$
For $n=5$ our method also proves easily    $d^{n-4}\rho_\varepsilon (a_{i_1,
\varepsilon} ,...,a_{i_k,\varepsilon} ) \rightarrow 0,$ but for the proof
of \, $d^{n-3}{\nabla{\rho}}_{\varepsilon} (a_{i_1,\varepsilon} ,...,a_{i_k,
\varepsilon} )$ goes to $0$ when $\e\to 0$, we need a more careful estimates of the rests
in Propositions \ref{p:32} and \ref{p:33} below.
\end{rem}
To proceed further, we introduce some notation.\\
Let, for $p\in\N^*$ and $\eta >0$ given
\begin{eqnarray*}
V_\varepsilon (p,\eta)&=\left\{\right.& u\in {\Sigma ^{+}(A_\varepsilon )}\,
\mbox{ s.t } \exists \, x_1,..., x_p \in A_\varepsilon, \exists \,
\lambda_1,...,\lambda_p  >0
\, \, \mbox{with} \nonumber \\
 & & \bigg|\bigg|u -C(p) \sum_{i=1}^pP_\varepsilon\delta_{( x_i,\lambda_i)}
\bigg|\bigg| <\eta ,\, \forall i \,
\lambda_i >\frac{1}{\eta } , \lambda_i \, d(x_i ,
\partial A_\varepsilon )\geq \frac{1}{\eta}, \nonumber \\
& &  \forall i\neq j \quad \varepsilon_{ij} = ( \frac{
\lambda_i}{\lambda_j}  + \frac{\lambda_j}{\lambda_i}
+\lambda_i\lambda_j|x_i - x_j|^2)^{-\frac{n-4}{2}}<\eta\left.  \right\},
\end{eqnarray*}
\noindent
where $\Sigma ^+ (A_{\varepsilon }) = \{ u \in E_\e /\, u > 0,\, \, ||u|| = 1\}$.\\
If a function $u$ belongs to $V_\varepsilon(p,\eta )$, then, for
$\eta >0$ small enough, the minimization problem
\begin{eqnarray}\label{e:31}
\min_{
 \alpha_i ,\lambda_i >0,\  x_i\in A_\varepsilon
} \left|\left|u - \sum_{i=1}^p\alpha_iP_\varepsilon\delta_{(
 x_i,\lambda_i)} \right|\right|
\end{eqnarray}
\noindent
has a unique solution, up to permutation (the proof of this fact is similar, up to minor modifications, to the corresponding statement for Laplacian operator in \cite {BC}).\\
Therefore, for $\varepsilon >0$ sufficiently small, Section 2
implies that $u_\varepsilon $( solution of $P_\varepsilon$ ) can
be uniquely written as
\begin{eqnarray}\label{e:32}
 \widetilde{u}_\varepsilon = u_\varepsilon/(||
u_\varepsilon||) = \sum_{i=1}^p\alpha_{i,\varepsilon}
P_\varepsilon\delta_{(a_{i,\varepsilon} , \lambda_{i,\varepsilon}
)} + v_\varepsilon,
\end{eqnarray}
\noindent
where $v_\varepsilon$ satisfies the following conditions :
$$
(V_0) \, \,  \left( v_\varepsilon , P_\varepsilon
\delta_{(a_{i,\varepsilon} ,\lambda_{i,\varepsilon} )}
\right) = \left(  v_\varepsilon ,\partial P_\varepsilon
\delta_{(a_{i,\varepsilon},\lambda_{i,\varepsilon} )}/(\partial
\lambda_i)\right)
=\left( v_\varepsilon ,\partial P_\varepsilon \delta_{(a_{i,\varepsilon}
 ,\lambda_{i,\varepsilon} )}/(\partial a_i)\right) = 0
$$
and $\alpha_{i,\varepsilon}$ satisfies :
$$
\left( J(u_{ \varepsilon})\right)^{n/(n-4)} \alpha_{j,\e}^{8/n-4}
= 1 +o(1) ,\forall j
$$
\noindent
In order to simplify the notations, in the sequel, we write $\alpha_i,a_i,
 \lambda_i, \delta_i $, $P\delta_i$ and $\th_i$ instead of $\alpha_{i,\varepsilon}$,
$a_{i,\varepsilon}$ , $\lambda_{i,\varepsilon}$,
$\delta_{(a_{i,\varepsilon} ,\lambda_{i,\varepsilon})}$, $P\delta_{(a_{i,
\varepsilon} ,\lambda_{i,\varepsilon})}$ and $\th_{(a_{i,\varepsilon} ,\lambda_{i,\varepsilon})}$ respectively and we also write
$u_\varepsilon$ instead of  $\widetilde{u}_\varepsilon.$ \\
First of all, we deal with the $v_\e$-part of $u_\e$.
\begin{pro}\label{p:31}
Let $v_\varepsilon$ be defined by \eqref{e:32}. Then we have the following estimate
$$
||v_{\varepsilon}|| \leq C
\begin{cases}
\sum_i\frac{1}{(\lambda_id_i)^{n-4}}   +
\sum\varepsilon_{ij}\left( Log\varepsilon_{ij}^{-1}
\right)^{\frac{n-4}{n}}&  \mbox{  if  } n<12   \\
\sum_i\frac{log(\l_i d_i)}{(\lambda_id_i)^{\frac{n+4}{2}}}
 +\sum\varepsilon_{ij}^{\frac{n+4}{2(n-4)}}\left(
Log\varepsilon_{ij}^{-1} \right)^{\frac{n+4}{2n}} & \mbox{  if  }
n\geq 12.
\end{cases}
$$
\end{pro}
\begin{pf}
From \eqref{e:32}, we derive that
\begin{eqnarray*}
\Delta^2 v_\varepsilon &=&J(u_{ \varepsilon })^{\frac{n}{n-4}}
\left(
 \sum\alpha_iP\delta_i +v_\varepsilon    \right)^{
\frac{n+4}{n-4}} - \sum\alpha_i\delta_i^{
\frac{n+4}{n-4}}\nonumber\\
&=&J(u_{ \varepsilon })^{\frac{n}{n-4}}\left[ \left( \sum\alpha_iP
\delta _i    \right)^{\frac{n+4}{n-4}} +\frac{n+4}{n-4}(
\sum\alpha_iP\delta_i)^{\frac{8}{n-4}}v_\varepsilon +
O(|v_\varepsilon |^{\frac{n+4}{n-4}} ) \nonumber \right. \\
&+& \left. O \left(Sup(\sum\alpha_iP\delta_i , v_\varepsilon )^{
\frac{8}{n-4}-1}|v_\varepsilon |^2   \right) \right]
-\sum\alpha_i\delta_i^{\frac{n+4}{n-4}}.
\end{eqnarray*}
\noindent
Thus, since $J(u_\varepsilon)$ is bounded,
$$
|| v_\varepsilon ||^2 =  J(u_{ \varepsilon
})^{\frac{n}{n-4}}\left[<f,v_\e > +\frac {n+4}{n-4}
\int_{A_\varepsilon}( \sum\alpha_i P\delta_i
)^{\frac{8}{n-4}}v_\varepsilon ^{2} \right] +
 O( ||v_\varepsilon ||^{\mbox{inf}(3,2n/(n-4)}),
$$
where
\begin{eqnarray}\label{e:33}
 <f ,v>  = \int_{A_\varepsilon}(\sum\alpha_iP\delta_i)^{\frac{n+4}{n-4}}v.
\end{eqnarray}
\noindent
Then
$$
 Q(v_\varepsilon , v_\varepsilon ) =  J(u_{
\varepsilon })^{\frac{n}{n-4}}<f,v_{ \varepsilon } > +
O(||v_\varepsilon ||^{\mbox{inf}(3,{2n}/{(n-4))}} ),
$$
\noindent
where
\begin{align*}
 Q(v_\varepsilon & , v_\varepsilon ) :=
|| v_\varepsilon ||^2 - \frac{n+4}{n-4}
 J(u_\varepsilon)^{\frac{n}{n-4}}\int_{A_\varepsilon}
(\sum_i\alpha_iP\delta_i)^{\frac{8}{n-4}}v_\varepsilon^2\\
 & =|| v_\varepsilon ||^2 - \frac{n+4}{n-4}
 J(u_\varepsilon)^{\frac{n}{n-4}}\sum_i\alpha_i ^{\frac{8}{n-4}}
\int_{A_\varepsilon} P\delta_i ^{\frac{8}{n-4}}v_\varepsilon^2+
O\left(\sum_{j\ne i} \int _{P\d_j \leq P\d_i} P\delta_i
^{\frac{8}{n-4}-1}P\d_j v_\varepsilon^2\right)
\end{align*}
Observe that,  since $J(u_\varepsilon)^{{n}/{(n-4)}}\alpha _i^{
{8}/(n-4)} = 1+o(1)$ and $\int (\d_i \d_j)^{n/(n-4)}=o(1)$, then
$Q(v,v)$ is close to
$$
|| v ||^2 - \frac{n+4}{n-4}\sum_i\int_{A_\varepsilon}
P\delta_i^{\frac{8}{n-4}}v^2
$$
\noindent
and therefore $Q$ is a positive definite quadratic form on $v$ (see \cite{BE}).
\\
Thus there exists $\a_0 >0$ such that
$$
\alpha_0||v_\varepsilon ||^2\leq  Q(v_\varepsilon , v_\varepsilon
) + O( ||v_\varepsilon ||^{\mbox{inf}(3,{2n}/{n-4}} ) =
J(u_\varepsilon) ^{\frac{n}{n-4}}<f,v_\varepsilon> \leq
C|f|||v_\varepsilon ||.
$$
\noindent
Hence
$$
||v_\varepsilon || \leq C'|f|.
$$
\noindent Now we estimate $|f|$. We have
$$
<f,v_\varepsilon > =
\sum_{i=1}^p\alpha_i^{\frac{n+4}{n-4}}\int_{A_\varepsilon}
P\delta_i^{\frac{n+4}{n-4}}v_\varepsilon + O\left(\sum_{j\ne i}
\int_{P\delta_j \leq
P\delta_i}P\delta_i^{\frac{8}{n-4}}P\delta_j|v_\varepsilon |
\right).
$$
\noindent
 Observe that
\begin{eqnarray*}
\int_{A_\varepsilon}P\delta_i^{\frac{n+4}{n-4}}v_\varepsilon
&=&\int _{ A\varepsilon }v_{ \varepsilon } \delta _{i}^{ \frac {n+4}{n-4}}
+ O\left(\int_{A_\varepsilon}\delta_i^{\frac{8}{n-4}}\th_i
|v_\varepsilon | \right)\nonumber\\
&=&O\left(\int_{B(a_i ,d_i)}\delta_i^{\frac{8}{n-4}}\th_i
 |v_\varepsilon | \right) +O\left(\int_{\R^n\setminus B(a_i ,d_i)}\delta_i^{
\frac{n+4}{n-4}} |v_\varepsilon | \right),
\end{eqnarray*}
\noindent
where $d_i = d(a_i ,\partial A_\varepsilon )$.\\
Thus
\begin{eqnarray}\label{e:34}
\int_{A_\varepsilon}P\delta_i^{\frac{n+4}{n-4}}v_\varepsilon =
O\left(||v_\varepsilon || \left[
|\th_i |_{L^\infty}\left(
\int_{ B(a_i ,d_i)}\delta_i^{\frac{16n}{n^2-16}}\right)^{\frac{n+4}{2n}}
+ \left(\int_{ B^c(a_i ,d_i)}
\delta_i^{\frac{2n}{n-4}}\right)^{ \frac {n+4}{2n}}\right]\right).
\end{eqnarray}
\noindent
Notice that
$$
\int_{\R^n\setminus B(a_i ,d_i)}\delta_i^{\frac{2n}{n-4}}  =
O\left( \frac{1}{(\lambda_id_i )^ { n}} \right)
$$
 and, since $|\theta_i|_{L^\infty} =O(\l_i ^{(4-n)/2}d_i ^{4-n})$,
we have
\begin{eqnarray*}
|\th_i|_{L^\infty}\left(\int_{B}\delta_i^{\frac{16n}{n^2-16}}
\right)^{ \frac{n+4}{2n}}  = O \left( \frac{1}{(\lambda _{i}
d_i)^{\frac{n+4}{2}}} + (\mbox{if } n=12) \frac {Log(\lambda _i
d_i)}{(\lambda _i d_i)^8} +(\mbox{if}\,  n<12) \frac {1}{(\lambda
_i d_i)^{n-4}} \right)
\end{eqnarray*}
where $B= B(a_i ,d_i)$. Therefore
\begin{eqnarray}\label{e:36}
\int_{A_\varepsilon}P\delta_i^{\frac{n+4}{n-4}}v_\varepsilon =  O
\left( ||v_\varepsilon ||\left[\frac{1}{(\lambda _{i}
d_i)^{\frac{n+4}{2}}} + (\mbox{if } n=12) \frac {Log(\lambda _i
d_i)}{(\lambda _i d_i)^8} +(\mbox{if}\,  n<12) \frac {1}{(\lambda
_i d_i)^{n-4}} \right]\right)
\end{eqnarray}
\noindent
We also have
\begin{eqnarray}\label{e:37}
\int_{P\delta_j\leq P\delta_i}P\delta_i^{\frac{8}{n-4}}P\delta_j
|v_\varepsilon | \leq ||v_{ \varepsilon }||\left[ \int_{P\delta_j
\leq P\delta_i}\left(P\delta_i) ^{\frac{8}{n-4}}P\delta_j\right)^{
\frac{2n}{n+4}}   \right]^{\frac{n+4}{2n}}.
\end{eqnarray}
\noindent
If $n\geq 12,$ we have $\frac{2n}{n+4}\geq\frac{n}{n-4}$ and thus
\begin{eqnarray}\label{e:38}
\int_{P\delta_j\leq P\delta_i}(P\delta_i^{\frac{8}{n-4}}P\delta_j
)^{\frac{2n}{n+4}} \leq  \int(\delta_i \delta_j )^{\frac{n}{n-4}}
= O\left( \e_{ij}^{\frac{n}{n-4}}Log \e _{ij} ^{-1} \right).
\end{eqnarray}
\noindent If $n\leq 11,$ we have $1<\frac{8}{n-4}$, thus, using
Holder's inequality, we derive that
\begin{eqnarray}\label{e:39}
 \int (P\delta_i^{\frac{8}{n-4}}
P\delta_j)^{\frac{2n}{n+4}} \leq c \left(\int
(\delta_i\delta_j)^{\frac{n}{n-4}} \right)^{ \frac{2(n-4)}{n+4}}
\leq  c\varepsilon _{ij}^{\frac{2n}{n+4}} \left(Log\varepsilon
_{ij}^{-1} \right)^{\frac{2(n-4)}{n+4}}
\end{eqnarray}
\noindent
Using \eqref{e:36}, \eqref{e:37}, \eqref{e:38} and  \eqref{e:39},
we easily deduce our proposition.
 \end{pf}

Next we will give useful expansions of gradient of $J$.
\begin{pro}\label{p:32}
For $n\geq 6,$  we have the following expansion
\begin{eqnarray*}
\left( \nabla  J(u_\varepsilon ),\lambda_i{\frac{
\partial P\delta_i}{\partial\lambda_i}}\right)
 & = & 2J(u_\varepsilon )c_1\left[-\frac{(n-4)}{2}
\alpha_i\frac{H_\varepsilon (a_i,a_i)}{\lambda_i^{n-4}  }(1+o(1))\right.\\
 & - & \left. \sum_{j\neq i}\alpha _j \left(\lambda _i \frac{\partial
\varepsilon_{ij}}{\partial\lambda_i}  + \frac{(n-4)}{2}
\frac{H_\varepsilon (a_i,a_j)}{(\lambda_i\lambda_j)^{
\frac{n-4}{2}}  }
\right)(1+o(1)) +R \right],
\end{eqnarray*}
\noindent where $R$ satisfies $ R =
o\left(\sum_k\frac{1}{(\lambda_kd_k )^{n-3}} + \sum_{k\ne
r}\varepsilon_{kr}^{(n-3)/(n-4)} \right). $
\end{pro}
\begin{pf}
We have
\begin{eqnarray}\label{e:310}
\nabla  J(u_\varepsilon ) = 2J(u_\varepsilon )\left[u_\varepsilon
  - J(u_\varepsilon )^\frac{n}{n-4}\Delta^{-2}(u_{\varepsilon}
 ^{\frac{n+4}{n-4}})\right].
\end{eqnarray}
\noindent
Thus, setting $\var_i = \lambda_i(\partial
 P\delta_i/\partial\lambda_i)$ and using Proposition \ref{p:31}, we have
\begin{align}\label{e:311}
\left( \nabla J(u_\varepsilon ),\var_i  \right)
 = & 2J(u_\varepsilon ) \left[ \sum \alpha_j( P\delta_j ,
\var_i) -J(u_\varepsilon )^\frac{n}{n-4} \left[  \int \left ( \sum
\alpha _j P \delta _ {j} \right )^ { \frac {n+4}{n-4}}
 \var_i\right.\right.\notag\\
 & + \left. \left.\frac{n+4}{n-4}\int (\sum\alpha_jP\delta_j)^\frac{8}{n-4}
v_\varepsilon \var_i \right] \right]+ R .
\end{align}
 Notice that if $ n \geq 6 $, we have
\begin{eqnarray}\label{e:313}
\int\left(\sum_{j=1}^p\alpha _j
P\delta_j\right)^\frac{8}{n-4}\var_i v_\e =\int (\alpha_i
P\delta_i)^\frac{8}{n-4}\var_i v_\e +O\left(\sum_{k\ne
r}\int_{P\d_r \leq P\d_k} P\d _k ^{\frac{8}{n-4}}P\d _r|v_\e|
\right).
\end{eqnarray}
Furthermore
\begin{align}\label{e:314}
 \int_{P\d_r\leq P\d_k} & P\delta_r P\delta_k^{\frac{8}{n-4}}
|v_\varepsilon | \leq   \int_{P\delta_r\leq |v_\varepsilon
|}P\d_k^{\frac{8}{n-4}}|v_\e|^2 +\int_{P\delta_r\leq
P\d_k}P\d_r^2 P\d_k^{\frac{8}{n-4}}\notag \\
   & = O\left( ||v_\varepsilon ||^2
+(\mbox{ if }n\leq 7) \left(\int (\delta _r\delta
_k)^{\frac{n}{n-4}} \right)^{\frac{2(n-4)}{n}}+(\mbox{ if }n\geq
8)\int (\delta _r\delta _k)^{\frac{n}{n-4}} \right)\notag\\
 & =  O\left(||v_\e||^2+(\mbox{ if }n\leq
7)\e_{kr}^2(log\e_{kr}^{-1})^{\frac{2(n-4)}{n}}+(\mbox{ if }n\geq
8)\e_{kr}^{\frac{n}{n-4}}log\e_{kr}^{-1}\right),
\end{align}
\begin{eqnarray}\label{e:316} \int_{\R^n\setminus B(a_i
,d_i)}P\delta_i^\frac{8}{n-4} v_\varepsilon\lambda_i\frac{\partial
P\delta_i}{\partial \lambda_i} = O\left(||v_\varepsilon ||
\frac{1} {(\lambda_id_i)^\frac{n+4}{2}}\right),
\end{eqnarray}
\begin{align}\label{e:315}
 \int_{B(a_i ,d_i)}&P\delta_i^{\frac{8}{n-4}} v_\varepsilon
\lambda_i\frac{\partial P\delta_i}{\partial\lambda_i} =
 - \int_{B_i}P\delta_i^{\frac{8}{n-4}}v_\varepsilon
\lambda_i\frac{\partial(\th_i)}{\partial\lambda_i} +
\int_{B_i}P\delta_i^{\frac{8}{n-4}}v_\varepsilon
\lambda_i\frac{\partial \delta_i}{\partial\lambda_i} \notag \\
 & \leq \frac{||v_ \varepsilon||}{(\l _i d_i ^2)^{\frac{n-4}{2}}}
\left( \int\delta^{\frac{16n}{n^2-16}} \right)^{\frac{n+4}{2n}}
+\int_{B_i}\delta_i^{\frac{8}{n-4}}
v_\varepsilon\lambda_i\frac{\partial \delta_i}{\partial\lambda_i}
+O\left(  \int_{B_i}\delta_i^{\frac{8}{n-4}}v_\varepsilon
\th_i \right) \notag \\
 & =O\left(||v_\varepsilon || \right)\left(
\frac{1}{(\lambda_id_i)^{\frac{n+4}{2}}} + (if n=12)
 \frac{Log(\lambda_id_i)}{(\lambda_id_i)^8} +
(if n\leq 11)\frac{1}{(\lambda_id_i)^{n-4} } \right).
\end{align}
Now, we need to estimate
\begin{align}\label{e:312}
\int (\sum_{j=1}^p\alpha _j P\delta_j)^\frac{n+4}{n-4} & \var_i
=\sum_{j=1}^p \int(\alpha _j P\delta_j)^\frac{n+4}{n-4}\var_i +
\frac{n+4}{n-4}\sum_{j\neq i}\int (
\alpha _iP\delta _i)^\frac{8}{n-4}\alpha _jP\delta_j \var_i \nonumber \\
 & + \sum_{k\ne r}O\left( (\mbox{ if }n\leq
7)\e_{kr}^2(log\e_{kr}^{-1})^{\frac{2(n-4)}{n}}+(\mbox{ if }n\geq
8)\e_{kr}^{\frac{n}{n-4}}log\e_{kr}^{-1}\right)
\end{align}
Now we observe that a computation similar to the one
performed in \cite{B} and \cite{R} shows that
\begin{align}
\bigl(P\d_i,\l_i\frac{\partial P\d_i}{\partial\l_i} \bigr)= &
\frac{n-4}{ 2}c_1\frac{H(a_i,a_i)}{\l_i
^{n-4}}+O\biggl(\frac{1}{(\l_i d_i)^{n-2}}\biggr)\label{e:1}\\
\int P\d_i^{\frac{n+4}{ n-4}}\l_i\frac{\partial P\d_i}{\partial\l_i} =
& 2\bigl(P\d_i,\l_i\frac{\partial P\d_i}{\partial\l_i} \bigr)+
O\biggl(\frac{1}{(\l_i d_i)^{n-2}}\biggr)\label{e:2}
\end{align}
and for $i\ne j$, we have
\begin{align}
\bigl(P\d_j,\l_i\frac{\partial P\d_i}{\partial\l_i}\bigr) =
c_1\biggl(\l_i & \frac{\partial\e_{ij}}{\partial\l_i}+ \frac{n-4}{
2}\frac{H(a_i,a_j)}{(\l_i\l_j)^{\frac{n-4}{ 2}}}\biggr)+
O\biggl(\sum_{k=i,j}\frac{1}{(\l_kd_k)^{n-2}}+
\e_{ij}^{\frac{n-2}{
n-4}}\biggr)\label{e:3}\\
\int P\d_j^{\frac{n+4}{ n-4}}\l_i\frac{\partial P\d_i}{
\partial\l_i}= & \bigl(P\d_j,\l_i\frac{\partial P\d_i}{
\partial\l_i}\bigr)+(\mbox{ if }
n\geq 8) O\left( \e_{ij}^{\frac{n}{n-4}}log(\e_{ij}^{-1})
+\frac{log(\l_jd_j)}{(\l_jd_j)^n}\right)\notag\\
 & \quad + (\mbox{ if } n\leq 7)O\left(\frac{\e_{ij}
(log\e_{ij}^{-1})^{\frac{n-4}{n}}}{(\l_jd_j)^{n-4}} \right)\label{e:5}\\
\int P\d_j\l_i\frac{\partial}{\partial\l_i} (P
\d_i)^{\frac{n+4}{n-4}} = & \bigl(P\d_j,\l_i\frac{\partial
P\d_i}{\partial\l_i}\bigr) +(\mbox{ if } n\geq 8)
O\left(\e_{ij}^{\frac{n}{n-4}}log(\e_{ij}^{-1})
 +\frac{log(\l_id_i)}{(\l_id_i)^n}\right)\notag\\
 & \quad + (\mbox{ if } n\leq 7)O\left(\frac{\e_{ij}
(log\e_{ij}^{-1})^{\frac{n-4}{n}}}{(\l_id_i)^{n-4}} \right).\label{e:6}
\end{align}
Now, the estimates \eqref{e:311},..., \eqref{e:6},
 and the fact that
 $J(u_\varepsilon )^\frac{n}{n-4}\alpha_j^ {\frac {8}{n-4}}
 = 1+o(1),$ Proposition \ref{p:32} follows.
\end{pf}
\begin{pro}\label{p:33}
For $n\geq 6$,  we have the following expansion
\begin{eqnarray*}
\left(\nabla J(u_\varepsilon ),\frac{1}{\lambda_i}\frac{
\partial P\delta_i}{\partial a_i}  \right) =
J(u )c_1 \left [\frac{\alpha_i}{\lambda_i^{n-1}}\frac{
\partial H_\varepsilon(a_i ,a_i )}{\partial a_i}(1+o(1))
-2 \sum_{j\neq i}\alpha _j \left( \frac{1}{\lambda_i} \frac{\partial
\varepsilon_{ij}}{\partial a_i}\nonumber \right. \right.\\
\left. \left.  -  \frac{1}{\lambda_i(\lambda_i \lambda_j)^{\frac{n-4}{2}}}
\frac{ \partial H_\varepsilon(a_i ,a_j )}{\partial a_i} \right)(1+o(1))
 + O\left(\sum\lambda_j|a_i - a_j
|\varepsilon_{ij}^{\frac{n-1}{n-4}}\right) \right]+R,
\end{eqnarray*}
\noindent
where $R$ is defined in Proposition \ref{p:32}.
\end{pro}
\begin{pf}
As in the proof of Proposition \ref{p:32} , we obtain
\eqref{e:311} but with ${\lambda_i}^{-1}(\partial P\delta_i)/
(\partial a_i)$ instead of $ \lambda _i (\partial
 P\delta_i)/(\partial \lambda_i)$.
Now, using Proposition 2.1 of \cite{BH1}, we derive that
\begin{align*}
\bigl(P\d_i,\frac{1}{\l_i}\frac{\partial P\d_i}{\partial
a_i}\bigr)= & -\frac{1}{ 2}\frac{c_1}{\l_i ^{n-3}}\frac{\partial
H}{\partial
a_i}(a_i,a_i) +O\bigl(\frac{1}{(\l_i d_i)^{n-2}}\bigr)\\
\int P\d_i ^{\frac{n+4}{ n-4}}\frac{1}{\l_i}\frac{\partial P\d_i}{\partial
a_i} = & 2\bigl(P\d_i,\frac{1}{\l_i}\frac{\partial P\d_i}{\partial
a_i} \bigr)+O\bigl(\frac{1}{(\l_i d_i)^{n-2}}\bigr)
\end{align*}
and for $ i\ne j$, we have
\begin{align*}
\bigl(P\d_j,\frac{1}{ \l_i}\frac{\partial P\d_i}{\partial
a_i}\bigr)= &
c_1\left(\frac{1}{\l_i}\frac{\partial\e_{ij}}{\partial a_i}
-\frac{1}{(\l_i \l_j)^{\frac{n-4}{
2}}}\frac{1}{\l_i}\frac{\partial H}{\partial a_i}(a_i,
a_j)\right)\\
 & \quad +O\biggl(\sum_{k=i,j}\frac{1}{(\l_kd_k)^{n-2}}+\e_{ij}^{\frac{n-1}{
n-4}}\l_j\mid a_i-a_j\mid\biggr)\\
\int P\d_j^{\frac{n+4}{ n-4}}\frac{1}{ \l_i}\frac{\partial P\d_i}{\partial
a_i}= & \bigl(P\d_j,\frac{1}{ \l_i}\frac{\partial P\d_i}{\partial
a_i}\bigr)+ (\mbox{ if } n\geq 8)
O\left(\e_{ij}^{\frac{n}{n-4}}log(\e_{ij}^{-1})
+\frac{log(\l_jd_j)}{(\l_jd_j)^n}\right)\\
 & \quad + (\mbox{ if } n\leq 7)O\left( \frac{\e_{ij}
(log\e_{ij}^{-1})^{\frac{n-4}{n}}}{(\l_jd_j)^{n-4}}
\right)\\
\int P\d_j\frac{1}{\l_i}\frac{\partial}{\partial a_i} (P\d_i)^{
\frac{n+4}{ n-4}}= & \bigl(P \d_j,\frac{1}{\l_i}\frac{\partial
P\d_i}{\partial a_i}\bigr)+(\mbox{ if } n\geq
8)O\left(\e_{ij}^{\frac{n}{n-4}}log(\e_{ij}^{-1})
+\frac{log(\l_id_i)}{(\l_id_i)^n}\right)\\
 & + (\mbox{ if } n\leq 7)O\left(\frac{\e_{ij}
(log\e_{ij}^{-1})^{\frac{n-4}{n}}}{(\l_id_i)^{n-4}} \right).
\end{align*}
Using the above estimates our proposition follows.
 \end{pf}\\
\noindent
Next we are going to give the proof of Theorem \ref{t:31}.\\
From Proposition \ref{p:32} we easily derive that $ p\geq 2$. Now for $i\in
\left\{1,...,p  \right\},$  we introduce the following condition
\begin{eqnarray}\label{e:317}
 {2^{-p-1}}\sum_{k\neq i}\varepsilon_{ki}\leq
\sum_{j=1}^p {H_\varepsilon(a_i ,a_j )}{(\lambda_i
\lambda_j)^\frac{4-n}{2}}.
\end{eqnarray}
\noindent
We divide the set $\{1,...,p\}$ into $T_1\cup T_2$ with
\begin{eqnarray*}
T_1 = \left\{ i / \, \, i \mbox{ satisfies } \eqref{e:317}
\right\}\quad \mbox{and}\quad T_2 = \left\{ i / \, \, i\mbox{ does
not satisfy }\eqref{e:317}\right\}.
\end{eqnarray*}
\noindent In $T_2$ we order the $\lambda_i 's$: $\lambda_{i_1}
\leq \lambda_{i_2} \leq ...\leq \lambda_{i_s}$.\\
We begin by proving the following Lemma:
\begin{lem}\label{l:34}
 For $ n\geq 6,$  we have the following estimate
$$
\sum_{j\in T_2,j\neq i}(\varepsilon_{ij} + {
(\lambda_id_i )^{4-n}} )  =R_1, \,\,\mbox{with}\,\,R_1 = o\left[\sum_{k\in T_1}\left(
 \frac{1}{(\lambda_kd_k)^{n-3}} +
 \sum_{r\neq k ,r\in T_1}\varepsilon_{kr}^\frac{n-3}
{n-4} \right)\right].\leqno{(a)}
$$
$$\sum_{k,r \in T_1} \e_{kr} = O\left(\sum _{k\in T_1}(\l_k
d_k)^{4-n}\right). \leqno{(b)}
$$
\end{lem}
\begin{pf}
We start by proving claim $(a)$. Using Proposition \ref{p:32}, we
derive that
\begin{align}\label{e:318}
0=\sum_{k=1}^s2^k\alpha_{i_k}(\nabla J(u_\varepsilon ),&
\lambda_{i_k}\frac{\partial P\delta_{i_k}}{\partial\l _{i_k}})
 = 2J( u_\varepsilon )c_1\sum_{k=1}^s \left [-\sum_{
j\neq i_k}2^k\alpha_j \alpha _{i_k} \lambda _{i_k}
\frac{\partial\varepsilon_{ji_k}}{
\partial\lambda_{i_k}}(1+o(1)) \right.\notag\\
 &-\left. \frac{n-4}{2}\sum_{j=1}^p2^k\alpha_j\alpha_{i_k}
\frac{H_\varepsilon(a_j ,a_{i_k})}{(\lambda_j
\lambda_{i_k})^\frac{n-4}{2}}(1+o(1)) +R \right].
\end{align}
\noindent
Notice that
\begin{eqnarray}\label{e:319}
 -\lambda_i\frac{\partial\varepsilon_{ij}}{\partial
\lambda_i} = \frac{n-4}{2}\varepsilon_{ij}\left (
1-\frac{2\lambda_j}{\lambda_i}\varepsilon_{ij}^\frac{2}{n-4}\right ).
\end{eqnarray}
\noindent Thus, if $\lambda_i\geq \lambda_j$ and $i$, $j$ in
$T_2$, we have
\begin{eqnarray}\label{e:320}
 -2\lambda_i\frac{\partial\varepsilon_{ij}}{\partial\lambda_i} -\lambda_j
\frac{\partial\varepsilon_{ij}}{\partial \lambda_j} \geq
-\lambda_i\frac{\partial\varepsilon_{ij}} {\partial\lambda_i}=
\frac{n-4}{2}\varepsilon_{ij}+O\left(\e_{ij}^{(n-2)/(n-4)}\right).
\end{eqnarray}
\noindent
For $j\in T_1$ and  $ i\in T_2,$  two cases may occur :\\
i/ $\frac{1}{2}d_j \leq d_i\leq2d_j$.\\
Using in this case the fact that $j$ satisfies \eqref{e:317} and $
H_\varepsilon( a_i ,a_k )\leq (d_id_k)^{\frac{-(n-4)}{2}}$, we
obtain
\begin{eqnarray}\label{e:321}
( \frac{\lambda_j}{\lambda_i} )^\frac{n-4}{2}\varepsilon_{ij}
\leq ( \frac{\lambda_j}{\lambda_i} )^\frac{n-4}{2}
c\sum_{k=1}^p\frac{H_\varepsilon(a_k ,a_j)}{
(\lambda_j\lambda_k)^\frac{n-4}{2}}\leq c
\sum_{k=1}^p\left((\lambda_id_i)(\lambda_kd_k)
 \right)^\frac{-(n-4)}{2}=o(1).
\end{eqnarray}
 \noindent
ii/ in other cases, we have $|a_i -a_j | \geq \frac{1}{2} \mbox
{max}(d_i ,d_j ),$ then
$$
\frac{\lambda_j}{\lambda_i}\varepsilon_{i j}^\frac{2}{n-4} \leq
\frac{\lambda_j}{ \lambda_i}(\lambda_i\lambda_j|a_i -a_j|^2)^{-1
}\leq(\lambda_i|a_i-a_j|)^{-2} = O\left({(\lambda_id_i)^{-2}}
\right)
$$
and \eqref{e:321} follows in this case.\\
Using \eqref{e:318}, \eqref{e:319}, \eqref{e:320}, and \eqref{e:321}, we see
that
$$
0\geq ((n-4)/2)\sum_{i\in T_2}\left( \sum_{j \neq
i}\varepsilon_{ij}(1+o(1))-2^p\sum_{j=1}^p {H_\varepsilon( a_i
,a_j)}{(\lambda_i\lambda_j)^\frac{4-n}{2}} + R \right).
$$
\noindent Since $i\in T_2$ and $H_\varepsilon(a_i,a_i) \sim{c}/{
d_i^{n-4}}$ for $\varepsilon$ small enough (this fact can be shown
as in the proof of Lemma 4.2 of \cite{AE}), then
$$
0\geq c\sum_{i\in T_2}\left( \sum_{j\neq i}\varepsilon_{ij}
+ (\lambda_id_i)^{4-n}\right)(1+o(1)) +R_1.
$$
 Therefore claim $(a)$ follows.\\
The proof of claim $(b)$ is based on \eqref{e:318} and claim
$(a)$.
 \end{pf}\\
\noindent Now, in $T_1$ we order all the $\lambda_id_i$'s:
 $\lambda_{j_1}d_{j_1}\leq \lambda_{j_2}d_{j_2} \leq ...
\leq\lambda_{j_q}d_{j_q}$. In order to simplify our notations, we
suppose that $T_1 = \{ 1, 2,...,q \}$ and $\lambda_1d_1\leq
\lambda_2d_2\leq ...\leq\lambda_qd_q$.  \\
Let us introduce the following sets:
\begin{align}
 K_0& = \left\{i\in T_1 / \exists k_1,...,k_m \in T_1 \, s.t.\,
\,  k_1 = i,...,k_m = 1 \mbox{ and } \frac{|a_{k_j} -a_{k_{j+1}}|}
{inf(d_{k_j} ,d_{k_{j+1}})}\leq C_0\right\}\label{e:322}\\
 B &= K_0\cap \left\{1,...,l \right\},\label{e:323}
\end{align}
\noindent where $l = max \{ i\in T_1\  s.t.\
{\lambda_kd_k}/{\lambda_{k-1}d_{k-1}} \leq C_1 \ \forall k\leq i
\}$ and $C_0 $ and $ C_1$ are positive constants chosen later.
\begin{lem}\label{l:35}
Let B be defined by \eqref{e:323}. Then, $\{1\} \varsubsetneq B$.
\end{lem}
\begin{pf}
We argue by contradiction.  We assume that $ B =\{ 1\}.$  \\
Using Proposition \ref{p:32}, and the fact that
$H_\varepsilon(a_i, a_i) \sim{c}/{d_i^{n-4}}$, we derive that
$$
0 = \left(\n J(u_\e ),\l _1\frac{\partial P\d _1}{\partial\l _1}\right) =
 2J(u)c_1 \left[ -\frac{(n-4)}{2}\alpha_1 \frac{
H_\varepsilon(a_1 ,a_1)}{\lambda_1^{n-4}}(1+o(1)) +
 O(\sum_{k\neq 1}\varepsilon_{k1})\right].
$$
\noindent
Thus
\begin{eqnarray}\label{e:324}
 0\leq -{c}{(\lambda_1d_1)^{4-n}} + O(
\sum_{k\neq 1}\varepsilon_{1k}).
\end{eqnarray}
\noindent
Observe that :\\
- for $k\in T_2,$ by Lemma \ref{l:34}, we have $
\varepsilon_{1k} = o\left((\l_1d_1)^{3-n})\right)$.\\
- For $k\in T_1,$  two cases may occur :\\
If $k>l$, then
\begin{align}\label{ee}
 \varepsilon_{1k} & \leq 2^{p+1}
\sum_{j=1}^p\frac{H_\varepsilon(a_k,a_j)}{
(\lambda_k\lambda_j)^\frac{n-4}{2}} \leq \sum_{j\in T_1} \frac{c}{
(\lambda_kd_k\lambda_jd_j)^\frac{n-4}{2}}+\sum_{j\in \{k\}\cup
T_2}\frac{c}{(\l_j d_j)^{n-4}}\notag\\
 & \leq C_1^\frac{4-n}{2}
((\lambda_ld_l)(\lambda_1d_1))^\frac{4-n}{2}+R_1.
\end{align}
Thus, using Lemma \ref{l:34} and the fact that $C_1$ large enough,
we obtain $\varepsilon_{1k} = o({(\lambda_1d_1
)^{4-n}})$.\\
Otherwise, we have $k\not\in K_0$, then $|a_1 -a_k|\geq C_0
\mbox{ inf } (d_1,d_k),$ then
$$
\varepsilon_{1k} \leq \left(\lambda_1 \lambda_k|a_1 - a_k|^2
\right)^\frac{4-n}{2} \leq C_0^{(4-n)/2}\left((
\lambda_1d_1)(\lambda_kd_k) \right)^\frac{4-n}{2}
 = o\left((\lambda_1d_1)^{4-n}   \right)
$$
\noindent
if we choose $C_0$ large enough .\\
Thus \eqref{e:324} yields a contradiction and our lemma follows.
\end{pf}\\
In order to finish the proof of  of Theorem \ref{t:31},  it is
sufficient to prove the following lemma.
\begin{lem}\label{l:36}
For $n \geq 6$, we have
$d^{n-4}\rho _B \longrightarrow 0$ and $d^{n-3}\nabla
 \rho _B \longrightarrow 0, \, $when $\varepsilon \longrightarrow 0,$\\
where $d = inf_{i \in B}d(a_i ,\partial A_\varepsilon )$
 and $\rho_B = \rho(a_{i_1},...,a_{i_m}),$ with \,
 $B = \{i_1,...,i_m \}$ the set defined by \eqref{e:323}.\\
\end{lem}
Before giving the proof of this lemma, we begin by studying the vector
$\Lambda $  defined by
\begin{eqnarray}\label{e:325}
 \Lambda  = \left({\lambda_{i_1}^\frac{4-n}{2}},
...,{\lambda_{i_m}^\frac{4-n}{2}}\right).
\end{eqnarray}
Let $ M_B = M( a_i ,i\in B )$  the matrix defined by \eqref{e:13}
and $\rho_B$ its least eigenvalue. We denote by  $e$ the
eigenvector associated to $\rho_B$. As in \cite{BLR}, we can
easily prove that all components of $e$ are strictly positive. Let
$\eta >0$ be such that for any $\gamma$ belongs to a neighborhood
$C(e,\eta )$ of $e$, we have
\begin{eqnarray}\label{e:326}
   ^T\gamma M_B\gamma  - \rho_B|\gamma |^2 \leq
 \frac{c_2}{d^{n-4}}|\gamma |^2 \mbox{ and } ^T
\gamma\frac{\partial M_B}{\partial a_i}\gamma =
 \left(\frac{\partial \rho _B}{\partial a_i} +
 o(\frac{1}{d^{n-3}})  \right)|\gamma |^2
 \end{eqnarray}
and for $\gamma\in (\R_+)^m\setminus C(e, \eta )$ , we have
\begin{eqnarray}\label{e:327}
 ^T\gamma M_B\gamma - \rho_B|\gamma |^2\geq
\frac{c_3|\gamma |^2}{d^{n-4}},
\end{eqnarray}
where
$$
C(e,\eta ) \subset \{ y \in (\R^*_+)^m s.t |
\frac{y}{|y|} - e| < \eta \}.
$$

\begin{lem}\label{l:37}
Let $\Lambda $  be defined by \eqref{e:325}.  Then $\Lambda \in C(e, \eta ).$
\end{lem}
\begin{pf}
We argue by contradiction. We assume that $\Lambda\in
(\R^*_+)^m\setminus C(e,\eta )$. Let
$$
\Lambda  (t) = |\Lambda |\frac{(1-t)\Lambda  +
 t|\Lambda  |e}{| (1-t)\Lambda + t|\Lambda |e |}
: =\frac{y(t)}{|y(t)|}.
$$
\noindent
From Proposition \ref{p:32}, we derive that
$$
(\n J(u_\e ),Z)_{|t=0} = -c\frac{d}{dt}
 \left(^T\Lambda  (t)M_B\Lambda (t)  \right)_{|t=0}  +
 O\left( \sum_{i\in B,j\in (T_1\setminus B)\cup T_2}
\varepsilon_{ij} \right)+R + o\left(  \frac{1}
{(\lambda_1d_1)^{n-4}} \right)
$$
where $Z$ is the vector field defined on the variables $\l $ along
the flow line defined by $\L (t)$.\\
Observe that
\begin{align*}
\frac{d}{dt} & \left( ^T\Lambda  (t)M_B\Lambda  (t) \right)
=\frac{d}{dt}\left(\frac{ ^T\Lambda  (t)M_B
\Lambda (t)}{|\Lambda (t) |^2} |\Lambda (0)|^2
 \right)\nonumber\\
&= |\Lambda (0)|^2\frac{d}{dt}\left( \rho_B
+\frac{(1-t)^2}{|y(t)|^2}(^T\Lambda  (0)M_B\Lambda (0)
 - \rho _B|\Lambda (0) |^2  ) \right)\nonumber\\
&= |\Lambda (0)|^2\left( \frac{2(1-t)}{|y(t)|^4}(^T
\Lambda (0)M_B\Lambda (0)  - \rho_B|\Lambda (0) |^2
 )(-(1-t)|\Lambda (0)|<e, \Lambda (0)> -t|\Lambda |^2) \right)
\end{align*}
\noindent
Thus
\begin{align*}
(\n J(u_\e ),Z)_{|t=0} =& -\frac{2c}{|\Lambda |^2}
(^T\Lambda M_B\Lambda   - \rho _B|\Lambda  |^2  )
(-|\Lambda |<e, \Lambda (0)>) \\
 &+ o\left( (\lambda_1d_1 )^{4-n}  \right) + O\left( \sum_{i\in B, j\in
(T_1\setminus B)\cup T_2}\varepsilon _{ij} \right).
\end{align*}
\noindent Since $|e| = 1,$  then there exists $k_0$ such that
$e_{k_0} \geq \frac{1}{m}.$  Thus
$$
<e,\Lambda (0)>=\sum_i e_i\Lambda _i\geq \frac{1}{m}\Lambda_{k_0}.
$$
\noindent
Using \eqref{e:327}, we obtain
\begin{eqnarray*}
(\n J(u_\e ),Z)_{|t=0}&\geq & cc_3d^{4-n}| \Lambda |\Lambda _{k_0}
+ o\left((\lambda_1d_1 )^{4-n}
 \right)  + O\left( \sum_{i\in B,j\in (T_1\setminus B)}
\varepsilon _{ij} \right)\nonumber \\
&\geq & c\left( \lambda _1 d_1 \lambda _{k_0} d_{k_0}
 \right)^{\frac {4-n}{2}} + o\left ( \frac {1}{\left(
\lambda _1 d_1 \right)^{n-4}}\right) + O \left( \sum _{i
\in B, \, j \in T_1 \setminus B} \varepsilon _{i,j} \right)
\end{eqnarray*}
\noindent
Observe that \\
- if $ j>l$,  since $ j \in T_1$, using \eqref{ee}, we have $
\varepsilon_{ij}=o (\lambda_1d_1\l_{k_0} d_{k_0} )^{(4-n)/2}$.\\
\noindent
 - if $j \not\in K_0$ and $j\leq l$
\begin{eqnarray*}
\varepsilon_{ij} \leq\left(\frac{1}{\lambda_i\lambda_j
|a_i-a_j|^2}\right)^\frac{n-4}{2} \leq \frac{
C_0^{(4-n)/2}}{(\lambda_1d_1\lambda_jd_j )^\frac{n-4}{2}} \leq
\frac{C_0^{\frac{4-n}{2}}C_1^{m-1}}{
(\lambda_1d_1\lambda_{k_0}d_{k_0} )^\frac{n-4}{2}} =o\left(
\frac{1}{(\lambda_1d_1\lambda_{k_0}d_{k_0} )^{\frac{n-4}{2}}}
\right)
\end{eqnarray*}
\noindent if we chose $C_0>>C_1$. Thus
$$
 0\geq \left( c(\lambda_1d_1\lambda_{k_0}d_{k_0}
)^\frac{4-n}{2}  \right) +o\left( (
\lambda_1d_1 )^{4-n}  \right) \geq \left(
(\lambda_1d_1 )^{4-n}  \right) \left(
 c(C_1^m )^\frac{4-n}{2} +o(1)   \right) >0.
$$
This yields a contradiction and our lemma follows.
 \end{pf}\\
\begin{pfn}{\bf Lemma \ref{l:36}}
Observe that, as in \eqref{e:324}, it is easy to prove that, for
$i,j$ in $T_1$, we have
$(\l_i/\l_j+\l_j/\l_i)\e_{ij}^{2/(n-4)}=o(1)$ and therefore
\begin{eqnarray}\label{*1}
\e_{ij}=\left(\l_i\l_j|a_i-a_j|^2\right)^{(4-n)/2}(1+o(1)).
\end{eqnarray}
Since $u_\e$ is a critical point of $J$, we have
\begin{eqnarray*}
  \sum_{i\in B}
\left(\n J(u_\e ),\l _i\frac{\partial P\d _i}{\partial \l
_i}\right)=0.
\end{eqnarray*}
Using Proposition \ref{p:32} and \eqref{*1}, we derive that
\begin{align*}
0 & = \sum_{i\in B}\left[\frac{ H_\varepsilon (a_i ,a_i
)}{\lambda_i^{n-4}}(1+ o(1)) - \sum_{j\ne i, j\in
B}(\varepsilon_{ij} -  \frac{ H_\varepsilon (a_i ,a_j
)}{(\lambda_i \lambda_j)^\frac{n-4}{2}})(1+ o(1))
+O(\sum_{j\in(T_1\setminus B)\cup T_2}
 \varepsilon_{ij} )+R \right]\\
 & =  ^T \Lambda M_B\Lambda +o\left( \frac{1}
{(\lambda_1d_1 )^{n-4}} \right) +R_1 + O(\sum_{j\in(T_1\setminus
B), i\in B} \varepsilon_{ij}).
\end{align*}
 Observe that, for $i\in B$ and $j\in T_1 \setminus B,$ we have,
as in the proof of Lemma \ref{l:35}, $\varepsilon_{ij}
=o({(\lambda_1d_1)^{4-n}})$. Thus
\begin{eqnarray}\label{e:328}
  0= ^T\Lambda M_B\Lambda + o\left( {(
\lambda_1d_1 )^{4-n}} \right).
\end{eqnarray}
We assume, arguing by contradiction, that $d^{n-4}\rho _B \not
\longrightarrow 0$, when $\varepsilon \longrightarrow 0.$
Therefore, there exists $C_4>0$ such that $|\rho _Bd^{n-4}|\geq C_4$ .\\
Now we distinguish two cases \\
{\bf $1^{st}$case :}$\rho _B>0$\\
In this case, we derive from \eqref{e:328}
$$
0\geq \rho _B|\Lambda |^2 + o\left(( \lambda_1d_1 )^{4-n} \right)
\geq C_4 |\Lambda |^2d^{4-n} +o\left(( \lambda_1d_1 )^{4-n}
\right)
>0.
$$
\noindent This yields a contradiction and we derive that
$d^{n-2}\rho_B
\longrightarrow 0$ in this case. \\
{\bf{$ 2^{nd}$ case :}}  $\rho_B \leq 0. $ In this case, we derive
from \eqref{e:326} and \eqref{e:328},
\begin{eqnarray*}
0\leq \rho _B|\Lambda |^2 +c_2 |
\Lambda |^2d^{4-n} +o\left((
\lambda_1d_1 )^{4-n} \right)& \leq&
 |\Lambda |^2d^{4-n}\left(
 \rho_B \, d^{n-4} +c_2   \right) +o\left(
(\lambda_1d_1 )^{4-n} \right)\nonumber \\
& \leq& |\Lambda |^2d^{4-n}(-C_4 +c_2) +o\left((\lambda_1d_1
)^{4-n} \right).
\end{eqnarray*}
\noindent If we choose $c_2 \leq \frac{1}{2}C_4,$ we obtain a
contradiction. Then $d^{n-4}\rho_B \longrightarrow 0,$  when
$\varepsilon
\longrightarrow 0$ also in this case. \\
Observe that, since $d^{n-4}\rho_B\longrightarrow 0,$  then there
exists $C_5>0$ such that $|a_i - a_j|\geq C_5 d,$   for any  $i ,j
\in B $ and  $i\neq j.$ \\
We assume, arguing by contradiction, that $d^{n-3}\nabla \rho_B \not
\longrightarrow 0$  when  $\varepsilon\longrightarrow 0.$\\
Since $u_\e$ is a critical point of $J$, we have $(\n
J(u_\e),\l^{-1}(\partial P\d_i)/(\partial a_i))=0$. For $i\in B$,
using Proposition \ref{p:33} and \eqref{*1}, we derive that
\begin{eqnarray*}
0 & = & ^T\Lambda \frac{\partial M_B}{\partial a_i} \Lambda +
O\left(\sum_{j\in T_1\setminus B} \left|\frac{
\partial\varepsilon_{ij}}{\partial a_i} -
 \frac{1}{( \lambda_i\lambda_j )^\frac{n-4}{2}}
\frac{\partial H_\varepsilon}{\partial a_i}(a_i ,a_j )
 \right|\right)\\
 & + & o\left( \frac{1}{d_i}\frac{1}
{(\lambda_1d_1 )^{n-4}}\right) +\lambda_iR
+ O\left( \sum_{j\in T_1} \lambda_i\lambda_j|a_i-a_j|
\varepsilon_{ij}^\frac{n-1}{n-4}\right) + O\left(
\sum_{j\in T_2}\lambda_i\varepsilon_{ij}\right).
\end{eqnarray*}
\noindent
Observe that:\\
- for $j\in T_2,$  we have, by Lemma \ref{l:34},
$\lambda_i\varepsilon_{ij} = o(d^{-1}(\l_1d_1)^{4-n}),$\\
- for $j\in T_1\setminus K_0$, since $(\partial/\partial
a_i)H(a_i, a_j) \leq d_i ^{-1}H(a_i,a_j)$, we have
\begin{eqnarray*}
|\frac{\partial\varepsilon_{ij}}{\partial a_i}| +
 \frac{1}{( \lambda_i\lambda_j )^\frac{n-4}{2}}
|\frac{\partial H_\varepsilon}{\partial a_i}(a_i,a_j)|
 &\leq& \frac{c}{(\lambda_i\lambda_j )^\frac{n-4}{2}}
 \left(  \frac{1}{|a_i-a_j|^{n-3}} +\frac{1}{
d_i|a_i-a_j|^{n-4}} \right)\\
&\leq& c C_0^{(4-n)/2}\left( ( \lambda_id_i\lambda_jd_j
)^\frac{4-n}{2}d_i
 \right)\\
&\leq& cC_0^{(4-n)/2}d^{-1} (\lambda_1d_1 )^{4-n} = o
\left(d^{-1}(\lambda_1d_1)^{4-n} \right).
\end{eqnarray*}
- for $j>l$, two cases may occur:\\
{\bf i/} The first case is when $d_j/2\leq d_i \leq 2d_j$. Using
the fact that $j\in T_1$, we obtain
\begin{align*}
\bigg|\frac{\partial\varepsilon_{ij}}{\partial a_i}\bigg| & +
 \frac{1}{( \lambda_i\lambda_j )^{\frac{n-4}{2}}}\bigg|\frac{
\partial H_\varepsilon}{\partial a_i}(a_i,a_j)\bigg|
\leq c\sqrt{\l_i \l_j}\e_{ij}^{\frac{n-3}{n-4}}
+\frac{c}{d_i}\frac{1}{(\l_id_i\l_jd_j)^{\frac{n-4}{2}}}\\
 & \leq  c\sqrt{\l_i\l_j}\sum_{k=1}^p\frac{1}
{(\l_kd_k\l_jd_j)^{\frac{n-3}{2}}}+ \frac{c}{C_1^\frac{n-4}{2}}
\left( \frac{1}{(\lambda_id_i)^{n-4}d_i}  \right)\\
 & \leq \frac{c\sqrt{\l_i\l_j}}{(\l_id_i\l_jd_j)^{1/2}}\sum_{k=1}^p\frac{1}
{(\l_kd_k)^{\frac{n-3}{2}(\l_jd_j)^{\frac{n-5}{2}}}}+
\frac{c}{C_1^\frac{n-4}{2}} \left( \frac{1}{(\lambda_1d_1)^{n-4}d}
\right).
\end{align*}
As in \eqref{ee} we obtain
$$\bigg|\frac{\partial\varepsilon_{ij}}{\partial a_i}\bigg| +
 \frac{1}{( \lambda_i\lambda_j )^{\frac{n-4}{2}}}\bigg|\frac{
\partial H_\varepsilon}{\partial a_i}(a_i,a_j)\bigg|=o\left(\frac{1}{d
(\l_1d_1)^{n-4}}\right).
$$
{\bf ii/} In other cases, we have $|a_i-a_j|\geq \max(d_i,d_j)/2$.
Then
\begin{align*}
\bigg|\frac{\partial\varepsilon_{ij}}{\partial a_i}\bigg| & +
 \frac{1}{( \lambda_i\lambda_j )^{\frac{n-4}{2}}}\bigg|\frac{
\partial H_\varepsilon}{\partial a_i}(a_i,a_j)\bigg|\leq
\frac{c}{(\l_i \l_j)^{(n-4)/2}|a_i-a_j|^{n-3}}+\frac{c}{d_i
(\l_id_i\l_jd_j)^{(n-4)/2}}\\
 & \leq \frac{c}{d_i (\l_id_i\l_jd_j)^{\frac{n-4}{2}}}\leq
\frac{c}{C_1^\frac{n-4}{2}} \left(
\frac{1}{(\lambda_1d_1)^{n-4}d}\right)=o\left(\frac{1}{d
(\l_1d_1)^{n-4}}\right).
\end{align*}
- for $j\in T_1$, as in the previous case, it is easy to see that
$$
\lambda_i\lambda_j|a_i-a_j|\varepsilon_{ij}^{
\frac {n-1}{n-4}} =  o\left({d(
\lambda_1d_1)^{3-n}}  \right),
$$
 Therefore, by \eqref{e:326}, we have
$$
0= ^T\Lambda \frac{\partial M_B}{\partial a_i}\Lambda
+  o\left(\frac{1}{d(\lambda_1d_1)^{n-4}}  \right)
= \left( \frac{\partial \rho_B}{\partial a_i} d^{n-3}
 +o(1) \right)\frac{|\Lambda |^2}{d^{n-3}} +
o\left(\frac{1}{d(\lambda_1d_1)^{n-4}}  \right).
$$
\noindent
Thus
$$
0\geq \left( |\n \rho|d^{n-3} +o(1) \right)\frac{|\Lambda |^2}
{d^{n-3}} + o\left(\frac{1}{d(\lambda_1d_1 )^{n-4}}  \right) \geq
C_6\frac{|\Lambda |^2}{d^{n-3}} + o\left(\frac{1}{d(
\lambda_1d_1)^{n-4}}  \right) >0
$$
\noindent
This yields a contradiction and our lemma follows.
\end{pfn}
\section{Proof of Theorem \ref{t:11} }
\mbox{}
Let us start by proving the following result :
\begin{thm}\label{t:41}
For $n\geq 5$,  let $C_0>0$ and let $(x_1 ,x_2 ,...,x_k)\in
A_\varepsilon^{k}$ such that
$$
d^{n-4}\rho_\varepsilon (x_1 ,...,x_k) \rightarrow 0\,\,
\mbox{  when }\,\, \varepsilon \rightarrow 0 \,\,\mbox{
and  }\,\,{|x_i - x_j|}\leq C_0 d , \forall  i , j,
$$
where $d = min\{d(x_i ,\partial A_\varepsilon ) / 1\leq i\leq k\}.$
\noindent
Then we have
$$
d^{n-3}\nabla{\rho_{\varepsilon}} (x_1,...,x_k) \not\longrightarrow 0,\, \mbox{
when } \varepsilon \rightarrow 0.
$$
\end{thm}
 \noindent
To prove Theorem \ref{t:41}, we introduce some notation and recall some result.\\
Let $(x_1,...,x_k)\in A_\varepsilon ^k $ such that
\begin{eqnarray}\label{e:41}
 d^{n-4}\rho _\varepsilon(x_1,...,x_k)
 \longrightarrow 0\mbox{   when  }\varepsilon
\longrightarrow 0\mbox{   and  }
{|x_i-x_j|}\leq C_0 d,\forall i,j,
\end{eqnarray}
\noindent
where $C_0$ is a fixed positive constant and $d= min_{1\leq i \leq k}
d(x_i ,\partial A_\varepsilon).$\\
We may assume, without loss of generality, that $d_1 =
inf_{1\leq i\leq k}d_i.$\\
Now we introduce the map
$$
{A_\varepsilon}\longrightarrow\widetilde{A}_\varepsilon, \quad
x\longmapsto \widetilde{x} =d_1^{-1} (x-x_1).
$$
As in the proof of Lemma 3.1 of \cite{AE}, we easily derive that
\begin{eqnarray}\label{e:42}
 \rho_\varepsilon(x_1,...,x_k) = d_{1}^{4-n}
\widetilde{\rho} _\varepsilon (0, \widetilde{x}_2 , . . . .,\widetilde{x}_k),
\end{eqnarray}
where $\widetilde{\rho}_\varepsilon$ is the function defined,
replacing $A_\varepsilon^k$ by $\widetilde{A}_\varepsilon^k$
in \eqref{e:13}, and $\widetilde{A_\varepsilon}$ converges in
the $C^1$ -topology on every compact set to $\Omega$ , where
$\Omega$ is a half-space or a strip .\\
Observe that $|\widetilde{x}_i| \leq C_0 ,\forall i\in
 \left\{ 2,...,k\right\}.$\\
Now we have the following Lemmas :
\begin{lem}\label{l:41}
For $\varepsilon >0,$  let
$$
F_k(\varepsilon )= \left\{ ( X_1,...,X_k )
\in\widetilde{A_\varepsilon^k} / \exists i\neq j \, \, s.t\, \,  X_i =X_j
\right\}
$$
\noindent
Then $\widetilde{\rho_\varepsilon}$ converges in the
$C^1$-topology to $\rho _\Omega$ , when$ \, \varepsilon
\longrightarrow 0$ , on every compact set which does not
intersect  V , where V is any neighborhood of
$F_k(\varepsilon )$ and $\rho _\Omega$ is the function
defined, replacing $ A_\varepsilon^k$ by $\Omega^k$ in \eqref{e:13}.
\end{lem}
\noindent
The proof of Lemma \ref{l:41} is similar to that of Lemma 4.1 in \cite{AE}.
\begin{lem}\label{l:42}
let $\rho _\Omega$ the function defined replacing
$ A_\varepsilon^k$ by $\Omega^k$ in (3). Then the map
$$
]0,1] \longrightarrow \R , \quad
t\longmapsto t^{n-4}\rho_\Omega( 0, tX_2,...,tX_k )
$$
\noindent
decreases when t decreases for any $X_2 , . . . .X_k \in\Omega.$
\end{lem}
\noindent
The same arguments in the proof of Lemma 4.5 in \cite{AE} prove easily our
lemma.\\
\begin{pfn}{Theorems \ref{t:41} and \ref{t:11}}
From \eqref{e:41}, \eqref{e:42} and Lemmas \ref{l:41} and \ref{l:42}, we
easily deduce  Theorem \ref{t:41}. Lastly,
 Theorem \ref{t:11} is an easy consequence of Theorems
\ref{t:21}, \ref{t:31} and \ref{t:41}.
\end{pfn}

\end{document}